\documentclass[a4paper,12pt]{article}
    \usepackage[top=2.5cm,bottom=2.5cm,left=2.5cm,right=2.5cm]{geometry}
    \usepackage{amsmath, amssymb}
    \allowdisplaybreaks

\usepackage{mathrsfs}

\usepackage{amsmath, epsfig, cite}
\usepackage{amssymb}
\usepackage{amsfonts}
\usepackage{latexsym}
\usepackage{graphicx}

\newtheorem{thm}{Theorem}[section]

\newtheorem{lem}[thm]{Lemma}

\newtheorem{con}[thm]{Conjecture}

\newtheorem{pro}[thm]{Proposition}

\newcommand{\qed}{{\hfill\rule{4pt}{7pt}}}
\def\pf{\noindent {\it Proof.} }

\numberwithin{equation}{section}

\makeatletter \@addtoreset{equation}{section} \makeatother

\begin{document}
\thispagestyle{empty}

\vspace*{20mm} \pagestyle{empty}

\begin{center}
{\Large\bf On a Conjecture about Tricyclic Graphs\\[3mm] with Maximal Energy\footnote{Supported
by NSFC.}}
\end{center}

\begin{center}
\small Xueliang Li,  Yongtang Shi\footnote{Corresponding author}, Meiqin Wei\\
\small Center for Combinatorics and LPMC-TJKLC\\
\small Nankai University, Tianjin 300071, China.\\
\small E-mail: lxl@nankai.edu.cn; shi@nankai.edu.cn;
weimeiqin@mail.nankai.edu.cn\\[2mm]
\small Jing Li\\
\small Department of Applied Mathematics \\
\small Northwestern Polytechnical University, Xi'an,
Shaanxi 710072, China \\
\small Email: jingli@nwpu.edu.cn

 (Received January 2, 2014)

\end{center}

\begin{center}
{\bf Abstract}
\end{center}

{\small For a given simple graph $G$, the energy of $G$, denoted by
$\mathcal {E}(G)$, is defined as the sum of the absolute values of
all eigenvalues of its adjacency matrix, which was defined by I.
Gutman. The problem on determining the maximal energy tends to be
complicated for a given class of graphs. There are many approaches
on the maximal energy of trees, unicyclic graphs and bicyclic
graphs, respectively. Let $P^{6,6,6}_n$ denote the graph with $n\geq
20$ vertices obtained from three copies of $C_6$ and a path
$P_{n-18}$ by adding a single edge between each of two copies of
$C_6$ to one endpoint of the path and a single edge from the third
$C_6$ to the other endpoint of the $P_{n-18}$. Very recently,
Aouchiche et al. [M. Aouchiche, G. Caporossi, P. Hansen, Open
problems on graph eigenvalues studied with AutoGraphiX, {\it Europ.
J. Comput. Optim.} {\bf 1}(2013), 181--199] put forward the
following conjecture: Let $G$ be a tricyclic graphs on $n$ vertices
with $n=20$ or $n\geq22$, then $\mathcal{E}(G)\leq
\mathcal{E}(P_{n}^{6,6,6})$ with equality if and only if $G\cong
P_{n}^{6,6,6}$. Let $G(n;a,b,k)$ denote the set of all connected
bipartite tricyclic graphs on $n$ vertices with three
vertex-disjoint cycles $C_{a}$, $C_{b}$ and $C_{k}$, where $n\geq
20$. In this paper, we try to prove that the conjecture is true for
graphs in the class $G\in G(n;a,b,k)$, but as a consequence we can
only show that this is true for most of the graphs in the class
except for 9 families of such graphs.}

\baselineskip=0.30in

\section{Introduction}
Let $G$ be a graph of order $n$ and $A(G)$ be the adjacency matrix of $G$. The
characteristic polynomial of $A(G)$ is defined as
$$
\phi(G,\lambda)=det(\lambda
I-A(G))=\sum\limits_{i=0}^{n}a_{i}\lambda^{n-i},
$$
which is called the \emph{characteristic polynomial} of $G$. The $n$
roots of the equation $\phi(G,\lambda)=0$, denoted by
$\lambda_{1},\lambda_{2},\cdots,\lambda_{n}$, are the
\emph{eigenvalues} of $G$. Since $A(G)$ is symmetric, all
eigenvalues of $G$ are real. It is well-known \cite{C.D1980} that if
$G$ is a bipartite graph, then
$$
\phi(G,\lambda)=\sum\limits_{i=0}^{\lfloor\frac{n}{2}\rfloor}a_{2i}\lambda^{n-2i}
=\sum\limits_{i=0}^{\lfloor\frac{n}{2}\rfloor}(-1)^{i}b_{2i}\lambda^{n-2i},
$$
where $b_{2i}=(-1)^{i}a_{2i}$ and $b_{2i}\geq 0$ for all $i=1,\cdots,\lfloor\frac{n}{2}\rfloor$.

The energy of $G$, denoted by $\mathcal{E}(G)$, is defined as
$$\mathcal{E}(G)=\sum\limits^{n}_{i=1}|\lambda_{i}|,$$ which was proposed by Gutman
in 1977 \cite{gutman1977}. The following formula is also well-known
\begin{equation*}
\mathcal{E}(G)=\frac{1}{\pi}\int^{+\infty}_{-\infty}\frac{1}{x^2}
\log|x^{n}\phi(G,i/x)|\mathrm{d}x,
\end{equation*}where $i^2=-1$. Moreover,
it is known from \cite{C.D1980} that the above equality can be expressed as the
following explicit formula:
\begin{equation*}
\mathcal{E}(G)=\frac{1}{2\pi}\int^{+\infty}_{-\infty}\frac{1}{x^{2}}\log
\left[\left(\sum\limits^{\lfloor n/2\rfloor}_{i=0}(-1)^{i}
a_{2i}x^{2i}\right)^{2}+\left(\sum\limits^{\lfloor
n/2\rfloor}_{i=0}(-1)^{i}
a_{2i+1}x^{2i+1}\right)^{2}\right]\mathrm{d}x,
\end{equation*}
where $a_{1}, a_{2}, \ldots, a_{n}$ are the coefficients of
$\phi(G,\lambda)$. It is also known \cite{IO} that for a bipartite
graph $G$, $\mathcal{E}(G)$ can be also expressed as the Coulson
integral formula
\begin{equation*}
\mathcal{E}(G)=\frac{2}{\pi}\int_{0}^{+\infty}\frac{1}{x^{2}}\ln\left[1+\sum\limits_{i=0}^{\lfloor\frac{n}{2}\rfloor}b_{2i}x^{2i}\right]dx.
\end{equation*}
For two bipartite graphs $G_{1}$ and $G_{2}$, if $b_{2i}(G_{1})\leq
b_{2i}(G_{2})$ hold for all
$i=1,2,\cdots,\lfloor\frac{n}{2}\rfloor$, we say that $G_{1}\preceq
G_{2}$ or $G_{2}\succeq G_{1}$. Moreover, if $b_{2i}(G_{1})<
b_{2i}(G_{2})$ holds for some $i$, we write $G_{1}\prec G_{2}$ or
$G_{2}\succ G_{1}$. Thus, for two bipartite graphs $G_{1}$ and
$G_{2}$, we can define the following quasi-order relation,
$$
G_{1}\preceq G_{2}\Rightarrow \mathcal{E}(G_{1})\leq
\mathcal{E}(G_{2}),\ \ G_{1}\prec G_{2}\Rightarrow
\mathcal{E}(G_{1})< \mathcal{E}(G_{2}).
$$For more results about graph energy, we refer the readers to two
surveys \cite{gutman2001,gutman&li&zhang2009} and the book \cite{LSG}.

It is quite interesting to study the extremal values of the energy
among some given classes of graphs, and characterize the
corresponding extremal graphs. In the meantime, a large number of
results were obtained on the minimal energies for distinct classes
of graphs, such as acyclic conjugated graphs \cite{zhang&li1999,
Li&Li2009}, bipartite graphs \cite{Li&Zhang&Wang2009}, unicyclic graphs \cite{Hou2001, Li&Zhou2005}, bicyclic
graphs \cite{Hou}, tricyclic graphs
\cite{Li&Li&Zhou2008, SLIXLI1} and tetracyclic graphs
\cite{SLIXLI2}. However, the maximal energy problem seems much more
difficult than the minimal energy problem. The commonly used
comparison method is the so-called quasi-order method. When the
graphs are acyclic, bipartite or unicyclic, it is almost always
valid. Nevertheless, for general graphs, the quasi-order method is
invalid. For these quasi-order incomparable problems, we found an
efficient way to determine which one attains the extremal value of
the energy, see \cite{HJL1,HJL2,HJLS1,HJLS2,HLS1,HLS2,HLSW}.

\begin{figure}[h,t,b,p]
\begin{center}
\includegraphics[scale = 0.7]{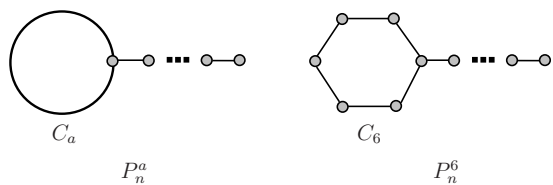}
\caption{Unicyclic graph $P_n^a$.}\label{fig1}
\end{center}
\end{figure}

Let $P_{n}$, $C_{n}$ and $S_n$ be a path, cycle and star garph with $n$ vertices,
respectively. Gutman \cite{gutman1977} first considered the extremal values of
energy of trees and showed that for any tree $T$ of order $n$, $\mathcal{E}(S_n)\leq
\mathcal{E}(T)\leq \mathcal{E}(P_n)$. Let $P_{n}^{a}$ be the graph obtained by
connecting a vertex of the cycle $C_{a}$ with a terminal vertex of the path
$P_{n-a}$ (as shown in Figure \ref{fig1}). In order to find lower and upper bounds
of the energy, Caporossi et al. \cite{CC} used the AGX system. They proposed a
conjecture on the maximal energy of unicyclic graphs.
\begin{con}
Among all unicyclic graphs on $n$ vertices, the cycle $C_n$ has maximal energy if $n\leq
7$ and $n=9,10,11,13$ and $15$\,. For all other values of $n$\,, the unicyclic graph
with maximal energy is $P_n^6$\,.
\end{con}

\begin{figure}[h,t,b,p]
\begin{center}
\includegraphics[scale = 0.7]{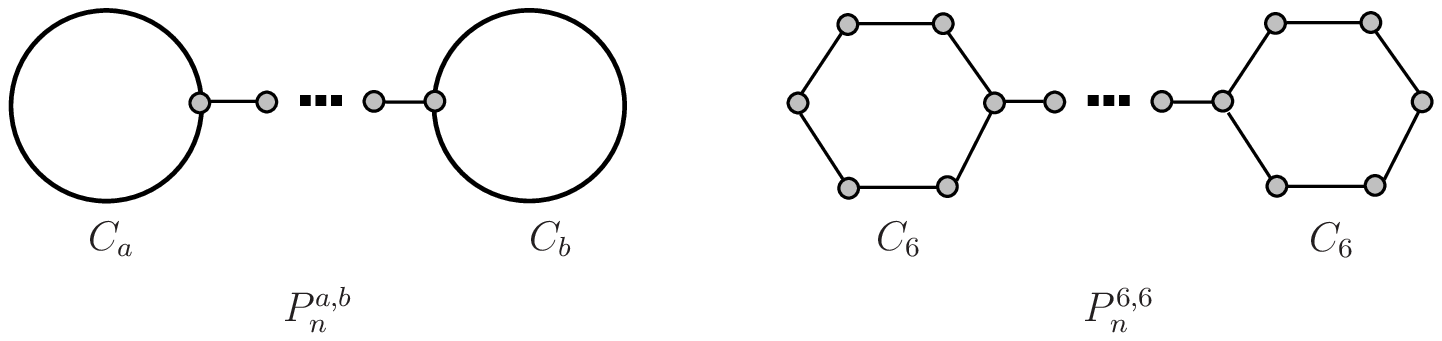}
\caption{Bicyclic graph $P_{n}^{a,b}$.}\label{fig2}
\end{center}
\end{figure}

In \cite{YIC}, Hou et al. proved a weaker result, namely that $\mathcal{E}(P_n^6)$
is maximal within the class of the unicyclic bipartite $n$-vertex graphs differing
from $C_n$\,. Huo et al. \cite{HLS1} and Andriantiana \cite{Andriantiana}
independently proved that $\mathcal{E}(C_n)<\mathcal{E}(P_n^6)$, and then completely
determined that $P_n^6$ is the only graph which attains the maximum value of the
energy among all the unicyclic bipartite graphs for $n=8,12,14$ and $n\geq 16$,
which partially solves the above conjecture. Finally, Huo et al. \cite{HLS2} and
Andriantiana and Wagner \cite{AW} completely solved this conjecture by proving the
following theorem, independently.
\begin{thm}
Among all unicyclic graphs on $n$ vertices, the cycle $C_n$ has maximal energy if $n\leq
7$ but $n\neq 4$, and $n=9,10,11,13$ and $15$\,;  $P_4^3$ has maximal energy if $n=4$\,.
For all other values of $n$\,, the unicyclic graph with maximal energy is $P_n^6$\,.
\end{thm}
The problem of finding bicyclic graphs with maximum energy was also
widely studied. Let $P_{n}^{a,b}$ (as shown in Figure \ref{fig2}) be
the graph obtained from cycles $C_{a}$ and $C_{b}$ by joining a path
of order $n-a-b+2$. Denote by $R_{a,b}$ the graph obtained from two
cycles $C_a$ and $C_b$ ($a, b\geq 10$ and $a \equiv b\equiv 2\,
(\,\textmd{mod}\, 4)$) by connecting them with an edge. In
\cite{gutman&vidovic2001}, Gutman and Vidovi\'{c} proposed a
conjecture on bicyclic graphs with maximal energy.

\begin{con}\label{conj}
For $n = 14$ and $n\geq 16$, the bicyclic molecular graph of order $n$ with maximal
energy is the molecular graph of the $\alpha,\beta$ diphenyl-polyene
$C_6H_5(CH)_{n-12}C_6H_5$, or denoted by $P^{6,6}_{n}$.
\end{con}

Furtula et al. \cite{FRG} showed by numerical computation that the
conjecture is true up to $n=50$. For bipartite bicyclic graphs, Li
and Zhang \cite{li&zhang2007} got the following result, giving a
partial solution to the above conjecture.
\begin{thm}\label{mainresultofliz}
If $G\in \mathscr{B}_{n}$, then $\mathcal{E}(G)\leq
\mathcal{E}(P^{6,6}_{n})$ with equality if and only if $G\cong
P^{6,6}_{n}$, where $\mathscr{B}_n$ denotes the class of all
bipartite bicyclic graphs but not the graph $R_{a, b}$.
\end{thm}
However, they could not compare the energies of $P^{6,6}_{n}$ and $R_{a,b}$. Furtula
et al. in \cite{FRG} showed by numerical computation that $\mathcal{E}(P^{6,6}_{n})>
\mathcal{E}(R_{a,b})$, which implies that the conjecture is true for bipartite
bicyclic graphs. They only performed the computation up to $a+b=50$. It is evident
that a solid mathematical proof is still needed. Huo et al. \cite{HJLS2} completely
solved this problem. However, the conjecture is still open for non-bipartite
bicyclic graphs.

\begin{thm} Let $G$ be any
connected, bipartite bicyclic graph with $n\,(\,n\geq 12)$ vertices.
Then $\mathcal{E}(G)\leq \mathcal{E}(P^{6,6}_{n})$ with equality if
and only if $G\cong P^{6,6}_{n}$.
\end{thm}

Actually, Wagner \cite{W} showed that the maximum value of the graph energy within
the set of all graphs with cyclomatic number $k$ (which includes, for instance,
trees or unicyclic graphs as special cases) is at most $4n/\pi + c_k$ for some
constant $c_k$ that only depends on $k$. However, the corresponding extremal graphs
are not considered.

\begin{figure}[h,t,b,p]
\begin{center}
\includegraphics[scale = 0.6]{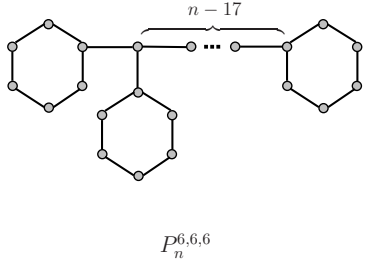}
\caption{Tricyclic graph $P_{n}^{6,6,6}$.}\label{fig12}
\end{center}
\end{figure}

The problem of finding the tricyclic graphs maximizing the energy
remains open. Gutman and Vidovi\'{c} \cite{gutman&vidovic2001}
listed some tricyclic molecular graphs that might have maximal
energy for $n\leq 20$. Very recently, in \cite{ACH}, experiments
using AutoGraphiX led us to conjecture the structure of tricyclic
graphs that presumably maximize energy for $n = 6,\ldots, 21$. For
$n\geq 22$, Aouchiche et al. \cite{ACH} proposed a general
conjecture obtained with AutoGraphiX. First, let $P^{6,6,6}_n$ (as
shown in Figure \ref{fig12}) denote the graph on $n\geq  20$
obtained from three copies of $C_6$ and a path $P_{n-18}$ by adding
a single edge between each of two copies of $C_6$ to one endpoint of
the path and a single edge from the third $C_6$ to the other
endpoint of the $P_{n-18}$.

\begin{con}\label{Con1.6}
Let $G$ be a tricyclic graphs on $n$ vertices with $n=20$ or $n\geq22$. Then
$\mathcal{E}(G)\leq \mathcal{E}(P^{6,6,6}_{n})$ with equality if and only if $G\cong
P^{6,6,6}_{n}$.
\end{con}

Let $G(n;a,b,k)$ denote the set of all connected bipartite tricyclic
graphs on $n$ vertices with three disjoint cycles $C_{a}$, $C_{b}$
and $C_{k}$, where $n\geq 20$. In this paper, we try to prove that
the conjecture is true for graphs in the class $G\in G(n;a,b,k)$,
but as a consequence we can only show that this is true for most of
the graphs in the class except for 9 families of such graphs.

\section{Preliminaries}
The following are the elementary results on the characteristic
polynomial of graphs and graph energy, which will be used in our
proof.
\begin{lem}\cite{C.D1980}\label{lemma2.1}
Let $uv$ be an edge of $G$. Then
$$
\phi(G,\lambda)=\phi(G-uv,\lambda)-\phi(G-u-v,\lambda)
-2\sum\limits_{C\in \varphi(uv)}\phi(G-C,\lambda),
$$
where $\varphi(uv)$ is the set of cycles containing $uv$. In
particular, if $uv$ is a pendant edge of $G$ with the pendant vertex
$v$, then
$$
\phi(G,\lambda)=\lambda\phi(G-v,\lambda)-\phi(G-u-v,\lambda).
$$
\end{lem}

\begin{lem}\label{lemma2.2}
Let $uv$ be an edge of a bipartite tricyclic graph $G$ which
contains three vertex-disjoint cycles. Then
$$
b_{2i}(G)=b_{2i}(G-uv)+b_{2i-2}(G-u-v)+2\sum\limits_{C_{l}\in
\varphi(uv)}(-1)^{1+\frac{l}{2}}b_{2i-l}(G-C_{l}),
$$
where $\varphi(uv)$ is the set of cycles containing $uv$. In
particular, if $uv$ is a pendant edge of $G$ with the pendant vertex
$v$, then
$$
b_{2i}(G)=b_{2i}(G-uv)+b_{2i-2}(G-u-v).
$$
\end{lem}

\pf By Lemma \ref{lemma2.1}, we have
$$
a_{2i}(G)=a_{2i}(G-uv)-a_{2i-2}(G-u-v)-2\sum\limits_{C_{l}\in
\varphi(uv)}a_{2i-l}(G-C_{l})
$$
and
\begin{eqnarray*}
(-1)^{i}a_{2i}(G)&=&(-1)^{i}a_{2i}(G-uv)+(-1)^{i-1}a_{2i-2}(G-u-v)\\[2mm]
&&+2\sum\limits_{C_{l}\in\varphi(uv)}(-1)^{1+\frac{l}{2}}
(-1)^{i-\frac{l}{2}}a_{2i-l}(G-C_{l}).\\
\end{eqnarray*}
\vskip -2em \noindent Since $b_{2i}=(-1)^{i}a_{2i}$, then the result
follows. \qed

From Sachs Theorem \cite{C.D1980}, we can obtain the following
properties for bipartite graphs.
\begin{pro}\label{Pro2.3}
(1). If $G_{1}$ and $G_{2}$ are both bipartite graphs,
then $b_{2k}(G_{1}\cup G_{2})=\sum\limits_{i=0}^{k}b_{2i}(G_{1})\cdot b_{2k-2i}(G_{2})$. \\
(2). Let $G$ and $G+e$ both be bipartite graphs, where $e\notin
E(G)$ and $G+e$ denotes the graph obtained from $G$ by adding the
edge $e$ to it. If either the length of any cycle
containing $e$ equals $2$ (mod $4$) or $e$ is not contained in any cycle, then $G\preceq G+e$. \\
(3). If $G_{0}$, $G_{1}$, $G_{2}$ are all bipartite and $G_{1}\preceq G_{2}$, since
$b_{2i}(G_{0})\geq 0$ and $b_{2i}(G_{1})\geq b_{2i}(G_{2})$ for all positive integer
$i$, we have $G_{0}\cup G_{1}\preceq G_{0}\cup G_{2}$. Moreover, for bipartite
graphs $G_{i}$, $G'_{i}$, $i=1,2$, if $G_{i}$ has the same order as $G'_{i}$ and
$G_{i}\preceq G'_{i}$, then $G_{1}\cup G_{2}\preceq G'_{1}\cup G'_{2}$.
\end{pro}

\begin{lem}\cite{IO}\label{lemma2.4}
Let $n=4k,4k+1,4k+2\ or\ 4k+3$. Then
\begin{eqnarray*}
P_{n}&\succ& P_{2}\cup P_{n-2}\succ P_{4}\cup P_{n-4}\succ\cdots\succ P_{2k}\cup P_{n-2k}\succ P_{2k+1}\cup P_{n-2k-1}\\
&\succ& P_{2k-1}\cup P_{n-2k+1}\succ\cdots\succ P_{3}\cup
P_{n-3}\succ P_{1}\cup P_{n-1}.
\end{eqnarray*}
\end{lem}

From the definition of $G(n;a,b,k)$, we know that $a$, $b$ and $k$
are all even. We will divide $G(n;a,b,k)$ into two categories
$G_{\uppercase\expandafter{\romannumeral1}}(n;a,b,k;l_{1},l_{2};l_{c})$
and
$G_{\uppercase\expandafter{\romannumeral2}}(n;a,b,k;l_{1},l_{2},l_{3})$
in the following.

We say that $H$ is the \emph{central structure} of $G$ if $G$ can be
viewed as the graph obtained from $H$ by planting some trees on it.
The central structures of
$G_{\uppercase\expandafter{\romannumeral1}}(n;a,b,k;l_{1},l_{2};l_{c})$
and
$G_{\uppercase\expandafter{\romannumeral2}}(n;a,b,k;l_{1},l_{2},l_{3})$
are
$\Theta_{\uppercase\expandafter{\romannumeral1}}(n;a,b,k;l_{1},l_{2};l_{c})$
 and
$\Theta_{\uppercase\expandafter{\romannumeral2}}(n;a,b,k;l_{1},l_{2},l_{3})$,
respectively.

\begin{figure}[h,t,b,p]
\begin{center}
\includegraphics[scale = 0.6]{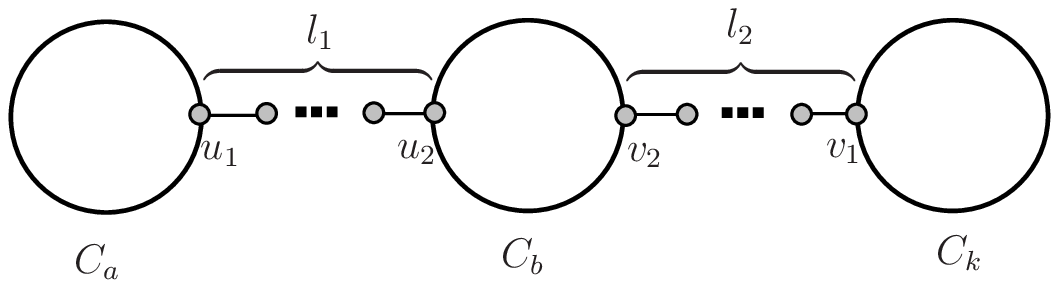}
\caption{$\Theta_{\uppercase\expandafter{\romannumeral1}}(n;a,b,k;l_{1},l_{2};l_{c})$.}
\label{fig3}
\end{center}
\end{figure}

$\Theta_{\uppercase\expandafter{\romannumeral1}}(n;a,b,k;l_{1},l_{2};l_{c})$
(as shown in Figure \ref{fig3}) is the set of all the elements of
$G(n;a,b,k)$ in which $C_{a}$ and $C_{b}$ are joined by a path
$P_{1}=u_{1}\cdots u_{2}$ ($u_{2}\in V(C_{b})$) with $l_{1}$
vertices,  $C_{k}$ and $C_{b}$ are joined by a path
$P_{2}=v_{1}\cdots v_{2}$ ($v_{2}\in V(C_{b})$) with $l_{2}$
vertices. In addition, the smaller part $u_{2}\cdots v_{2}$ of
$C_{b}$ has $l_{c}$ vertices. Note that when $u_2=v_2$, we have
$l_c=1$.

\begin{figure}[h,t,b,p]
\begin{center}
\includegraphics[scale = 0.6]{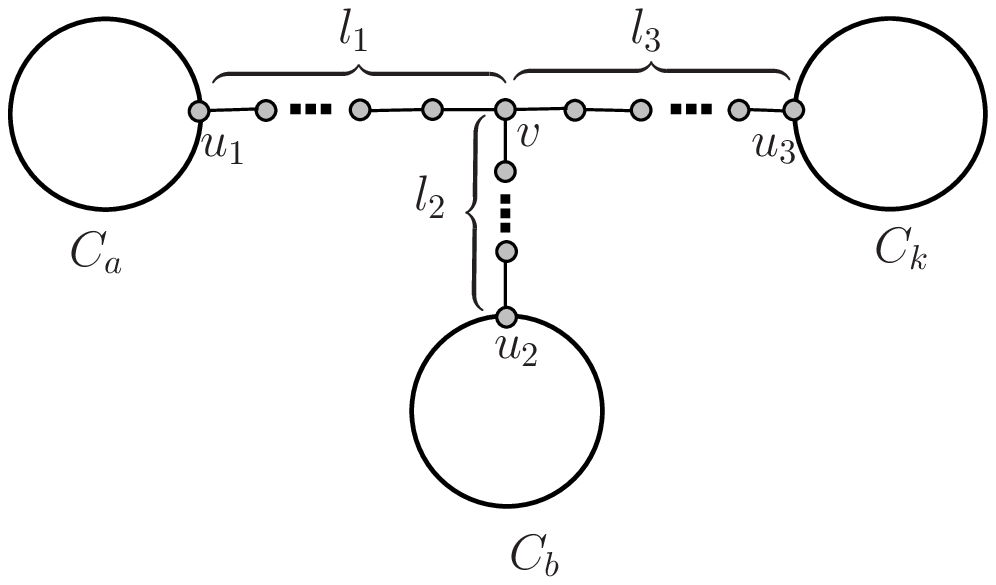}
\caption{$\Theta_{\uppercase\expandafter{\romannumeral2}}(n;a,b,k;l_{1},l_{2},l_{3})$.}
\label{fig4}
\end{center}
\end{figure}

$\Theta_{\uppercase\expandafter{\romannumeral2}}(n;a,b,k;l_{1},l_{2},l_{3})$
(as shown in Figure \ref{fig4}) is also a subset of $G(n;a,b,k)$.
For any $G\in
\Theta_{\uppercase\expandafter{\romannumeral2}}(n;a,b,k;l_{1},l_{2},l_{3})$,
$G$ has a center vertex $v$, $C_{a}$, $C_{b}$ and $C_{k}$ are joined
to $v$ by paths $P_{1}=u_{1}\cdots v$ ($u_{1}\in V(C_{a})$),
$P_{2}=u_{2}\cdots v$ ($u_{2}\in V(C_{b})$), $P_{3}=u_{3}\cdots v$
($u_{3}\in V(C_{k})$), respectively. The number of vertices of
$P_{1}$, $P_{2}$ and $P_{3}$ are $l_{1}$, $l_{2}$ and $l_{3}$,
respectively.

It is easy to verify that
\begin{eqnarray*}
G(n;a,b,k)&=&G_{\uppercase\expandafter{\romannumeral1}}(n;a,b,k;l_{1},l_{2};l_{c})
\cup\
G_{\uppercase\expandafter{\romannumeral2}}(n;a,b,k;l'_{1},l'_{2},l'_{3}).
\end{eqnarray*}

Now we define two special graph classes $\Gamma_{1}$ and
$\Gamma_{2}$ as follows.

 $\Gamma_{1}$ consists of graphs $G$ with the following four different
possible forms:\\
(\romannumeral 1)\ \ \ $G\in \Theta_{\uppercase\expandafter{\romannumeral 1}}(n;a,4,k;l_{1},l_{2};2)$,
where $a\geq 8$, $k\geq 8$, $2\leq l_{1}\leq 3$, $2\leq l_{2}\leq 3$. \\
(\romannumeral 2)\ \ $G\in
\Theta_{\uppercase\expandafter{\romannumeral 1}}(n;a,b,k;l_{1},l_{2};2)$,
where $a\geq 8$, $b\geq 6$, $k\geq 8$, $2\leq l_{1}\leq 3$, $2\leq l_{2}\leq 3$ and $l_{1}=l_{2}=3$
is not allowed. \\
(\romannumeral 3) $G\in \Theta_{\uppercase\expandafter{\romannumeral 1}}(n;4,b,k;l_{1},l_{2};2)$,
where $b\geq6$, $k\geq6$, $2\leq l_{1}\leq 3$ and $2\leq l_{2}\leq 3$. \\
(\romannumeral 4)\  $G\in
\Theta_{\uppercase\expandafter{\romannumeral
1}}(n;a,b,4;l_{1},l_{2};2)$, where $2\leq l_{2}\leq 3$.\\
Whereas $\Gamma_{2}$ consists of graphs $G$ with the following five different possible forms:\\
(\romannumeral 1)\ \ \ $G\in
\Theta_{\uppercase\expandafter{\romannumeral
2}}(n;a,b,k;2,l_{2},l_{3})$,
where $a\geq 8$. \\
(\romannumeral 2)\ \ $G\in \Theta_{\uppercase\expandafter{\romannumeral 2}}(n;a,b,k;3,3,3)$,
where $a\geq k\geq b\geq 8$. \\
(\romannumeral 3) $G\in \Theta_{\uppercase\expandafter{\romannumeral 2}}(n;a,4,k;l_{1},3,l_{3})$. \\
(\romannumeral 4) $G\in \Theta_{\uppercase\expandafter{\romannumeral 2}}(n;a,4,k;l_{1},2,l_{3})$. \\
(\romannumeral 5)\ \ $G\in
\Theta_{\uppercase\expandafter{\romannumeral 2}}(n;a,4,k;3,4,3)$,
where $a\geq k\geq 6$.

In this paper, we first try to find the graphs with maximal energy
among the two categories of $G(n;a,b,k)$:
$G_{\uppercase\expandafter{\romannumeral1}}(n;a,b,k;l_{1},l_{2};l_{c})$
and
$G_{\uppercase\expandafter{\romannumeral2}}(n;a,b,k;l_{1},l_{2},l_{3})$,
respectively. Then, we will obtain that
$P^{6,6,6}_{n}=\Theta_{\uppercase\expandafter{\romannumeral
2}}(n;6,6,6;n-17,2,2)$ has the maximal energy among all graphs in
$G(n;a,b,k)$ except for two classes $\Gamma_{1}$ and $\Gamma_{2}$.
Our main result is stated as follows, which gives support to
Conjecture \ref{Con1.6}.
\begin{thm}\label{mainthm}
For any tricyclic bipartite graph $G\in G(n;a,b,k)\setminus
(\Gamma_{1}\cup \Gamma_{2})$, $\mathcal {E}(G)\leq \mathcal
{E}(P^{6,6,6}_{n})$ and the equality holds if and only if $G\cong
P^{6,6,6}_{n}$.
\end{thm}

\section{Proof of Theorem \ref{mainthm}.}

By repeatedly applying the recursive formula of $b_{2i}(G)$ in Lemma
\ref{lemma2.2} and the third property in Proposition \ref{Pro2.3},
we obtain the following two lemmas.

\begin{lem}
If $G\in G_{\uppercase\expandafter{\romannumeral
1}}(n;a,b,k;l_{1},l_{2};l_{c})\setminus
\Theta_{\uppercase\expandafter {\romannumeral
1}}(n;a,b,k;l'_{1},l'_{2};l'_{c})$, then there exists a graph $G'\in
\Theta_{\uppercase\expandafter{\romannumeral
1}}(n;a,b,k;l'_{1},l'_{2};l'_{c})$ such that $G\prec G'$, i.e., the
graph with maximal energy among graphs in
$G_{\uppercase\expandafter{\romannumeral
1}}(n;a,b,k;l_{1},l_{2};l_{c})$ must belong to
$\Theta_{\uppercase\expandafter{\romannumeral
1}}(n;a,b,k;l'_{1},l'_{2};l'_{c})$.
\end{lem}

\begin{lem}
If $G\in G_{\uppercase\expandafter{\romannumeral
2}}(n;a,b,k;l_{1},l_{2},l_{3})\setminus
\Theta_{\uppercase\expandafter {\romannumeral
2}}(n;a,b,k;l'_{1},l'_{2},l'_{3})$, then there exists a graph $G'\in
\Theta_{\uppercase\expandafter{\romannumeral
2}}(n;a,b,k;l'_{1},l'_{2},l'_{3})$ such that $G\prec G'$, i.e., the
graph with maximal energy among graphs in
$G_{\uppercase\expandafter{\romannumeral
2}}(n;a,b,k;l_{1},l_{2},l_{3})$ must belong to
$\Theta_{\uppercase\expandafter{\romannumeral
2}}(n;a,b,k;l'_{1},l'_{2},l'_{3})$.
\end{lem}

From the results above, we know that the graph with maximal energy
among graphs in $G(n;a,b,k)$ must belong to
$\Theta_{\uppercase\expandafter{\romannumeral
1}}(n;a,b,k;l_{1},l_{2};l_{c})$ or
$\Theta_{\uppercase\expandafter{\romannumeral
2}}(n;a,b,k;l_{1},l_{2},l_{3})$. Therefore, in the following, we
will find the graph with maximal energy among graphs in
$\Theta_{\uppercase\expandafter{\romannumeral
1}}(n;a,b,k;l_{1},l_{2};l_{c})$ and
$\Theta_{\uppercase\expandafter{\romannumeral
2}}(n;a,b,k;l_{1},l_{2},l_{3})$.

\begin{lem}\label{lemma3.3}
For any graph $G\in \Theta_{\uppercase\expandafter{\romannumeral
1}}(n;a,b,k;l_{1},l_{2};l_{c})$, there exists a graph $H\in
\Theta_{\uppercase\expandafter{\romannumeral 1}}(n;a,\\b,k;l_{1},l_{2};2)$ such that
$G\preceq H$.
\end{lem}

\pf We distinguish the following two cases:

\textbf{Case 1.} $l_c=1$.

\begin{figure}[h,t,b,p]
\begin{center}
\includegraphics[scale = 0.6]{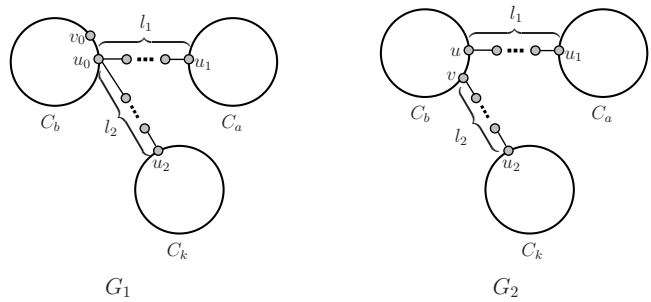}
\caption{Graphs for Lemma \ref{lemma3.3}.} \label{fig5}
\end{center}
\end{figure}

For fixed parameters $n$, $a$, $b$, $k$, $l_{1}$ and $l_{2}$, let
$G_{1}\in \Theta_{\uppercase\expandafter{\romannumeral
1}}(n;a,b,k;l_{1},l_{2};1)$ and
$G_{2}=\Theta_{\uppercase\expandafter{\romannumeral
1}}(n;a,b,k;l_{1},l_{2};2)$ (as shown in Figure \ref{fig5}). It
suffices to show that $G_{1}\preceq G_{2}$.

By Lemma \ref{lemma2.2} we have
\begin{eqnarray*}
b_{2i}(G_{1})&=&b_{2i}(G_{1}-u_{0}v_{0})+b_{2i-2}(G_{1}-u_{0}-v_{0})
+(-1)^{1+\frac{b}{2}}2b_{2i-b}(G_{1}-C_{b})\\
&=&b_{2i}(G_{1}-u_{0}v_{0})+b_{2i-2}(P_{a+l_{1}-2}^{a}\cup P_{k+l_{2}-2}^{k}\cup P_{b-2})\\
&&+(-1)^{1+\frac{b}{2}}2b_{2i-b}(P_{a+l_{1}-2}^{a}\cup
P_{k+l_{2}-2}^{k})
\end{eqnarray*}
and
\begin{eqnarray*}
b_{2i}(G_{2})&=&b_{2i}(G_{2}-uv)+b_{2i-2}(G_{2}-u-v)
+(-1)^{1+\frac{b}{2}}2b_{2i-b}(G_{2}-C_{b})\\
&=&b_{2i}(G_{2}-uv)+b_{2i-2}(P_{a+l_{1}-2}^{a}\cup P_{k+l_{2}-2}^{k}\cup P_{b-2})\\
&&+(-1)^{1+\frac{b}{2}}2b_{2i-b}(P_{a+l_{1}-2}^{a}\cup
P_{k+l_{2}-2}^{k}).
\end{eqnarray*}
Therefore, it suffices to show that $b_{2i}(G_{1}-u_{0}v_{0})\leq
b_{2i}(G_{2}-uv)$. By Lemma \ref{lemma2.2} we have
\begin{eqnarray*}
b_{2i}(G_{1}-u_{0}v_{0})&=&b_{2i}(P_{a+l_{1}-2}^{a}\cup P_{k+b+l_{2}-2}^{k})+b_{2i-2}(P_{a+l_{1}-3}^{a}\cup P_{k+l_{2}-2}^{k}\cup P_{b-1})\\
b_{2i}(G_{2}-uv)&=&b_{2i}(P_{a+l_{1}-2}^{a}\cup P_{k+b+l_{2}-2}^{k})+b_{2i-2}(P_{a+l_{1}-3}^{a}\cup P_{k+b+l_{2}-3}^{k})\\
&=&b_{2i}(P_{a+l_{1}-2}^{a}\cup P_{k+b+l_{2}-2}^{k})
+b_{2i-2}(P_{a+l_{1}-3}^{a}\cup P_{k+l_{2}-2}^{k}\cup P_{b-1})\\
&&+b_{2i-4}(P_{a+l_{1}-3}^{a}\cup P_{k+l_{2}-3}^{k}\cup P_{b-2}).
\end{eqnarray*}
Since $b_{2i-4}(P_{a+l_{1}-3}^{a}\cup P_{k+l_{2}-3}^{k}\cup
P_{b-2})\geq 0$, then we obtain $b_{2i}(G_{1}-u_{0}v_{0})\leq
b_{2i}(G_{2}-uv)$.

\textbf{Case 2.} $l_c\geq 2$.

\begin{figure}[h,t,b,p]
\begin{center}
\includegraphics[scale = 0.7]{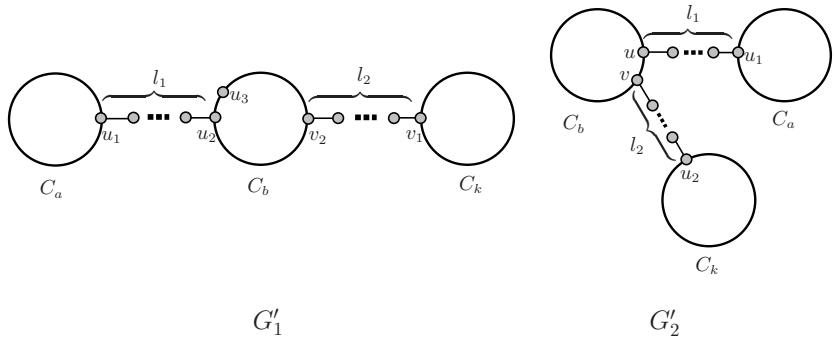}
\caption{Graphs for Lemma \ref{lemma3.3}.} \label{fig6}
\end{center}
\end{figure}

For fixed parameters $n$, $a$, $b$, $k$, $l_{1}$ and $l_{2}$, let
$G'_{1}\in \Theta_{\uppercase\expandafter{\romannumeral
1}}(n;a,b,k;l_{1},l_{2};l_{c})$ and
$G'_{2}\in\Theta_{\uppercase\expandafter{\romannumeral
1}}(n;a,b,k;l_{1},l_{2};2)$ (as shown in Figure \ref{fig6}, where
$u_3$ belongs to the part of $C_b$ with length $b-l_c+1$). It
suffices to show that $G'_{1}\preceq G'_{2}$.
 \\
By Lemma \ref{lemma2.2} we have\\
$b_{2i}(G'_{1})=b_{2i}(G'_{1}-u_{2}u_{3})+b_{2i-2}(G'_{1}-u_{2}-u_{3})
+(-1)^{1+\frac{b}{2}}2b_{2i-b}(G'_{1}-C_{b});$

\noindent$b_{2i}(G'_{2})=b_{2i}(G'_{2}-uv)+b_{2i-2}(G'_{2}-u-v)
+(-1)^{1+\frac{b}{2}}2b_{2i-b}(G'_{2}-C_{b}).$\\
Since
$(-1)^{1+\frac{b}{2}}2b_{2i-b}(G'_{1}-C_{b})=(-1)^{1+\frac{b}{2}}2b_{2i-b}(G'_{2}-C_{b})$,
we only need to compare
$b_{2i}(G'_{1}-u_{2}u_{3})+b_{2i-2}(G'_{1}-u_{2}-u_{3})$ with
$b_{2i}(G'_{2}-uv)+b_{2i-2}(G'_{2}-u-v)$. By applying Lemma
\ref{lemma2.2} repeatedly, we have
\begin{eqnarray*}
&&b_{2i}(G'_{1}-u_{2}u_{3})+b_{2i-2}(G'_{1}-u_{2}-u_{3})\\
&=&b_{2i}(P_{a+l_{1}+l_{c}-3}^{a}\cup P_{b+k+l_{2}-l_{c}-1}^{k}) +b_{2i-2}(P_{a+l_{1}+l_{c}-4}^{a}\cup P_{k+l_{2}-2}^{k}\cup P_{b-l_{c}} )\\
&&+b_{2i-2}(P_{a+l_{1}-2}^{a}\cup P_{k+l_{2}-2}^{k}\cup
P_{b-2})+b_{2i-4}(P_{a+l_{1}-2}^{a}\cup P_{k+l_{2}-3}^{k}\cup P_{l_{c}-2}\cup
P_{b-l_{c}-1}),
\end{eqnarray*}
and
\begin{eqnarray*}
&&b_{2i}(G'_{2}-uv)+b_{2i-2}(G'_{2}-u-v)\\
&=&b_{2i}(P_{a+l_{1}+l_{c}-3}^{a}\cup P_{b+k+l_{2}-l_{c}-1}^{k}) +b_{2i-2}(P_{a+l_{1}+l_{c}-4}^{a}\cup P_{b+k+l_{2}-l_{c}-2}^{k})\\
&&+b_{2i-2}(P_{a+l_{1}-2}^{a}\cup P_{k+l_{2}-2}^{k}\cup P_{b-2})\\
&=&b_{2i}(P_{a+l_{1}+l_{c}-3}^{a}\cup P_{b+k+l_{2}-l_{c}-1}^{k}) +b_{2i-2}(P_{a+l_{1}+l_{c}-4}^{a}\cup P_{k+l_{2}-2}^{k}\cup P_{b-l_{c}})\\
&&+b_{2i-4}(P_{a+l_{1}+l_{c}-4}^{a}\cup P_{k+l_{2}-3}^{k}\cup P_{b-l_{c}-1})+b_{2i-2}(P_{a+l_{1}-2}^{a}\cup P_{k+l_{2}-2}^{k}\cup P_{b-2})\\
&=&b_{2i}(P_{a+l_{1}+l_{c}-3}^{a}\cup P_{b+k+l_{2}-l_{c}-1}^{k}) +b_{2i-2}(P_{a+l_{1}+l_{c}-4}^{a}\cup P_{k+l_{2}-2}^{k}\cup P_{b-l_{c}})\\
&&+b_{2i-2}(P_{a+l_{1}-2}^{a}\cup P_{k+l_{2}-2}^{k}\cup P_{b-2})+b_{2i-4}(P_{a+l_{1}-2}^{a}\cup P_{k+l_{2}-3}^{k}\cup P_{l_{c}-2}\cup P_{b-l_{c}-1})\\
&&+b_{2i-6}(P_{a+l_{1}-3}^{a}\cup P_{k+l_{2}-3}^{k}\cup
P_{l_{c}-3}\cup P_{b-l_{c}-1}).
\end{eqnarray*}
Since $b_{2i-6}(P_{a+l_{1}-3}^{a}\cup P_{k+l_{2}-3}^{k}\cup
P_{l_{c}-3}\cup P_{b-l_{c}-1})\geq 0$, we have
$b_{2i}(G'_{1}-u_{2}u_{3})+b_{2i-2}(G'_{1}-u_{2}-u_{3})\leq
b_{2i}(G'_{2}-uv)+b_{2i-2}(G'_{2}-u-v)$.

Thus, the proof is complete. \qed

\begin{thm}\label{Thm3.4}
For any graph $G\in \Theta_{\uppercase\expandafter{\romannumeral
2}}(n;a,b,k;l_{1},l_{2},l_{3})\setminus \Gamma_{2}$, there exists a
graph $H\in \Theta_{\uppercase\expandafter{\romannumeral
2}}(n;6,6,6;l'_{1},l'_{2},l'_{3})$ such that $G\preceq H$.
\end{thm}

\pf Without loss of generality, we may assume that $a\geq k\geq b$.
It is obvious that $l_1, l_2, l_3\geq 2$. We distinguish the
following cases:

\textbf{Case 1.} $\begin{cases}
l_1+a-2\geq 7\\
l_2+b-3\geq 7\\
l_3+k-2\geq 7
\end{cases}$

In this case, considering the values of $l_1$, $l_2$ and $l_3$, we
distinguish the following four subcases.

\textbf{Subcase 1.1.} $l_1\geq3$, $l_2\geq4$ and$l_3\geq3$.

\begin{figure}[h,t,b,p]
\begin{center}
\includegraphics[scale = 0.7]{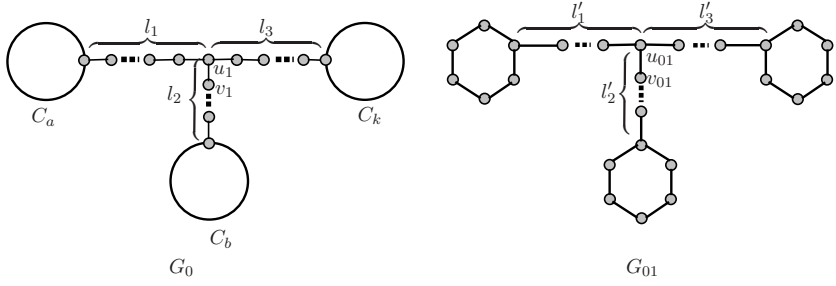}
\caption{Graphs for Subcase 1.1.} \label{fig7}
\end{center}
\end{figure}

For any values of $l_1$, $l_2$ and $l_3$, let $G_1\in
\Theta_{\uppercase\expandafter{\romannumeral
2}}(n;a,b,k;l_{1},l_{2},l_{3})$ and $G_{01}\in
\Theta_{\uppercase\expandafter{\romannumeral
2}}(n;6,6,\\6;l'_{1},l'_{2},l'_{3})$ (as shown in Figure
\ref{fig7}), where $l'_{1}=a+l_1-6$, $l'_{2}=b+l_2-6$ and
$l'_{3}=k+l_3-6$.

By Lemma \ref{lemma2.2}, we have
\begin{eqnarray*}
b_{2i}(G_{1})&=&b_{2i}(G_{1}-u_{1}v_{1})+b_{2i-2}(G_{1}-u_{1}-v_{1})\\
&=&b_{2i}(P^{a,k}_{a+k+l_{1}+l_{3}-3}\cup P^{b}_{b+l_{2}-2})+b_{2i-2}(P^{a}_{a+l_{1}-2}\cup P^k_{k+l_{3}-2}\cup P^b_{b+l_{2}-3}),\\
b_{2i}(G_{01})&=&b_{2i}(G_{01}-u_{01}v_{01})+b_{2i-2}(G_{01}-u_{01}-v_{01})\\
&=&b_{2i}(P^{6,6}_{l'_{1}+l'_{3}+9}\cup P^{6}_{l'_{2}+4})+b_{2i-2}(P^{6}_{l'_{1}+4}\cup P^6_{l'_{3}+4}\cup P^6_{l'_{2}+3})\\
&=&b_{2i}(P^{6,6}_{a+k+l_{1}+l_{3}-3}\cup
P^{6}_{b+l_{2}-2})+b_{2i-2}(P^{6}_{a+l_{1}-2}\cup
P^6_{k+l_{3}-2}\cup P^6_{b+l_{2}-3}).
\end{eqnarray*}
By Proposition \ref{Pro2.3}, we can obtain that $G_1\preceq G_{01}$.

\textbf{Subcase 1.2.} $l_1=2$, $l_2\geq4$, $l_3\geq3$ or $l_1\geq3$,
$l_2\geq4$, $l_3=2$.

The graphs in this subcase belong to $\Gamma_2(\romannumeral 1)$, so we do not consider them.

\textbf{Subcase 1.3.} $l_1\geq3$, $l_2=3$, $l_3\geq3$.

It is easy to verify that $b\geq 8$ and then we have $a\geq k\geq
b\geq 8$. Let $G_2\in\Theta_{\uppercase\expandafter{\romannumeral
2}}(n;a,b,k;l_{1},l_{2},l_{3})$ and $G_{02}\in
\Theta_{\uppercase\expandafter{\romannumeral
2}}(n;6,6,6;l'_{1},l'_{2},l'_{3})$, where
$l'_{1}=a+l_1-6$, $l'_{2}=b+l_2-6$ and $l'_{3}=k+l_3-6$.
If $l_1>3$ or $l_3>3$, then with similar analysis in Subcase 1.1, we
have
\begin{eqnarray*}
b_{2i}(G_{2})&=&b_{2i}(P^{k,b}_{k+b+l_{2}+l_{3}-3}\cup P^{a}_{a+l_{1}-2})+b_{2i-2}(P^{a}_{a+l_{1}-3}\cup P^k_{k+l_{3}-2}\cup P^b_{b+l_{2}-2}),\\
b_{2i}(G_{02})&=&b_{2i}(P^{6,6}_{l'_{2}+l'_{3}+9}\cup P^{6}_{l'_{1}+4})+b_{2i-2}(P^{6}_{l'_{1}+3}\cup P^6_{l'_{3}+4}\cup P^6_{l'_{2}+4})\\
&=&b_{2i}(P^{6,6}_{k+b+l_{2}+l_{3}-3}\cup
P^{6}_{a+l_{1}-2})+b_{2i-2}(P^{6}_{a+l_{1}-3}\cup
P^6_{k+l_{3}-2}\cup P^6_{b+l_{2}-2}).
\end{eqnarray*}
By Proposition \ref{Pro2.3}, we can obtain that $G_2\preceq G_{02}$.

If $l_1=l_2=l_3=3$, the graphs in this case belong to $\Gamma_2(\romannumeral 2)$, so we do not consider them.

\textbf{Subcase 1.4.} $l_1\geq3$, $l_2=2$, $l_3\geq3$ or
$l_1=l_3=2$, $l_2\geq4$ or $l_1=2$, $2\leq l_2\leq3$, $l_3\geq3$ or
$l_1\geq3$, $2\leq l_2\leq3$, $l_3=2$ or $l_1=l_3=2$, $2\leq
l_2\leq3$.

The graphs in this case belong to $\Gamma_2(\romannumeral 1)$, so we do not consider them.\\

\textbf{Case 2.} $\begin{cases}
l_1+a-2\leq 6\\
l_2+b-3\geq 7\\
l_3+k-2\geq 7
\end{cases}$

In this case, it is easy to verify that $a\leq6$, from which we have
$b\leq k \leq a \leq 6$. If $a=b=k=6$, it follows that this lemma
holds. Hence, we consider the following subcases.

\textbf{Subcase 2.1.} $a=k=6$, $b=4$.

It is easy to verify that $l_2\geq6$ and $l_3\geq3$. For any values
of $l_2$ and $l_3$, let $G_3\in
\Theta_{\uppercase\expandafter{\romannumeral
2}}(n;6,4,6;2,l_{2},l_{3})$ and $G_{03}\in
\Theta_{\uppercase\expandafter{\romannumeral
2}}(n;6,6,6;2,l'_{2},l'_{3})$, where $l'_2=l_2-2$ and
$l'_3=l_3$. By Lemma \ref{lemma2.2}, we have
\begin{eqnarray*}
b_{2i}(G_{3})&=&b_{2i}(P^{6,6}_{l_{3}+11}\cup P^{4}_{l_{2}+2})+b_{2i-2}(C_6\cup P^{6}_{l_{3}+4}\cup P^4_{l_{2}+1}),\\
b_{2i}(G_{03})&=&b_{2i}(P^{6,6}_{l_{3}+11}\cup
P^{6}_{l_{2}+2})+b_{2i-2}(C_6\cup P^{6}_{l_{3}+4}\cup
P^6_{l_{2}+1}).
\end{eqnarray*}
By Proposition \ref{Pro2.3}, we can obtain that $G_3\preceq G_{03}$.

\textbf{Subcase 2.2.} $a=6$, $k=b=4$.

It is easy to verify that $l_2\geq6$ and $l_3\geq5$. For any values
of $l_2$ and $l_3$, let $G_3\in
\Theta_{\uppercase\expandafter{\romannumeral
2}}(n;6,4,4;2,l_{2},l_{3})$ and $G_{03}\in
\Theta_{\uppercase\expandafter{\romannumeral
2}}(n;6,6,6;2,l'_{2},l'_{3})$, where $l'_2=l_2-2$ and
$l'_3=l_3-2$. By Lemma \ref{lemma2.2}, we have
\begin{eqnarray*}
b_{2i}(G_{4})&=&b_{2i}(P^{6,4}_{l_{3}+9}\cup P^{4}_{l_{2}+2})+b_{2i-2}(C_6\cup P^{4}_{l_{3}+2}\cup P^4_{l_{2}+1}),\\
b_{2i}(G_{04})&=&b_{2i}(P^{6,6}_{l_{3}+9}\cup
P^{6}_{l_{2}+2})+b_{2i-2}(C_6\cup P^{6}_{l_{3}+2}\cup
P^6_{l_{2}+1}).
\end{eqnarray*}
By Proposition \ref{Pro2.3}, we can obtain that $G_4\preceq G_{04}$.

\textbf{Subcase 2.3.} $a=k=b=4$.

It is easy to verify that $l_1\leq 4$, $l_2\geq 6$ and $l_3\geq 5$.
If $l_1=4$, let $G_5\in \Theta_{\uppercase\expandafter{\romannumeral
2}}(n;4,4,4;4,l_{2},l_{3})$ and $G_{05}\in
\Theta_{\uppercase\expandafter{\romannumeral
2}}(n;6,6,6;2,l'_{2},l'_{3})$, where $l'_2=l_2-2$ and
$l'_3=l_3-2$. By Lemma \ref{lemma2.2}, we have
\begin{eqnarray*}
b_{2i}(G_{5})&=&b_{2i}(P^{4,4}_{l_{3}+9}\cup
P^{4}_{l_{2}+2})+b_{2i-2}(P^4_6\cup P^{4}_{l_{3}+2}\cup
P^4_{l_{2}+1}),\\
b_{2i}(G_{05})&=&b_{2i}(P^{6,6}_{l_{3}+9}\cup
P^{6}_{l_{2}+2})+b_{2i-2}(C_6\cup P^{6}_{l_{3}+2}\cup
P^6_{l_{2}+1}).
\end{eqnarray*}
Also, $\phi(P^4_6;\lambda)=\lambda^6-6\lambda^4+6\lambda^2$ and
$\phi(C_6;\lambda)=\lambda^6-6\lambda^4+9\lambda^2-4$. It follows
that $P^4_6\prec C_6$. By Proposition \ref{Pro2.3}, we can obtain
that $G_5\preceq G_{05}$.

If $l_1<4$, graphs in this case belong to $\Gamma_2(\romannumeral 3)$ or
$\Gamma_2(\romannumeral 4)$, so we do not consider them.

\textbf{Case 3.} $\begin{cases}
l_1+a-2\geq 7\\
l_2+b-3\geq 7\\
l_3+k-2\leq 6
\end{cases}$

In this case, it is easy to verify that $k\leq 6$. If $b\leq k\leq a\leq 6$, with
similar analysis in Case 2 we obtain that this lemma holds. Then we consider the
case $a>6\geq k\geq b$. Without considering graphs with forms
$\Gamma_2(\romannumeral 1)$, $\Gamma_2(\romannumeral 3)$ and $\Gamma_2(\romannumeral
4)$, there are only two subcases as follows.

\textbf{Subcase 3.1.} $a>6$, $k=6$, $b=6\ or\ 4$.

It is easy to verify that $l_2\geq4$ and $l_3=2$. We have $l_1\geq3$
since we do not consider graphs with form $\Gamma_2(\romannumeral
1)$. For any values of $l_1$ and $l_3$, let $G_6\in
\Theta_{\uppercase\expandafter{\romannumeral
2}}(n;a,b,6;l_{1},l_2,2)$ and $G_{06}\in
\Theta_{\uppercase\expandafter{\romannumeral
2}}(n;6,6,6;l'_{1},l'_{2},2)$, where $l'_1=a+l_1-6$ and
$l'_2=b+l_2-6$. By Lemma \ref{lemma2.2}, we have
\begin{eqnarray*}
b_{2i}(G_{6})&=&b_{2i}(P^{a,6}_{a+l_{1}+5}\cup
P^{b}_{b+l_{2}-2})+b_{2i-2}(P^{a}_{a+l_{1}-2}\cup C_6\cup
P^b_{b+l_{2}-3}),\\
b_{2i}(G_{06})&=&b_{2i}(P^{6,6}_{a+l_{1}+5}\cup
P^{6}_{b+l_{2}-2})+b_{2i-2}(P^{6}_{a+l_{1}-2}\cup C_6\cup
P^6_{b+l_{2}-3}).
\end{eqnarray*}
By Proposition \ref{Pro2.3}, we can obtain that $G_6\preceq G_{06}$.

\textbf{Subcase 3.2.} $a>6$, $k=b=4$.

It is easy to verify that $l_2\geq6$ and $l_3\leq4$. We have
$l_1\geq3$ since we do not consider graphs with form
$\Gamma_2(\romannumeral 1)$. For
any values of $l_1$ and $l_3$, let $G_7\in
\Theta_{\uppercase\expandafter{\romannumeral
2}}(n;a,4,4;l_{1},l_2,4)$ and $G_{07}\in
\Theta_{\uppercase\expandafter{\romannumeral
2}}(n;6,6,6;l'_{1},l'_{2},2)$, where $l'_1=a+l_1-6$ and
$l'_2=l_2-2$. By Lemma \ref{lemma2.2}, we have
\begin{eqnarray*}
b_{2i}(G_{7})&=&b_{2i}(P^{a,4}_{a+l_{1}+5}\cup P^{4}_{l_{2}+2})+b_{2i-2}(P^{a}_{a+l_{1}-2}\cup P^{4}_{6}\cup P^4_{l_{2}+1}),\\
b_{2i}(G_{07})&=&b_{2i}(P^{6,6}_{a+l_{1}+5}\cup
P^{6}_{l_{2}+2})+b_{2i-2}(P^{6}_{a+l_{1}-2}\cup C_{6}\cup
P^6_{l_{2}+1}).
\end{eqnarray*}
Also, $P^4_6\prec C_6$. By Proposition \ref{Pro2.3}, we can obtain
that $G_7\preceq G_{07}$.

\textbf{Case 4.} $\begin{cases}
l_1+a-2\geq 7\\
l_2+b-3\leq 6\\
l_3+k-2\geq 7
\end{cases}$

It is easy to verify that $b\leq6$. Without considering graphs with
forms $\Gamma_2(\romannumeral 1)$, $\Gamma_2(\romannumeral 3)$,
$\Gamma_2(\romannumeral 4)$ and $\Gamma_2(\romannumeral 5)$, we can
distinguish this case into the following four subcases.

\textbf{Subcase 4.1.} $b=6$, $l_2=3$, $l_1\geq3$ and $l_3\geq3$.

For any values of $l_1$ and $l_3$, let $G_8\in
\Theta_{\uppercase\expandafter{\romannumeral 2}}(n;a,6,k;l_{1},3,l_3)$ and $G_{08}\in
\Theta_{\uppercase\expandafter{\romannumeral 2}}(n;6,6,6;l'_{1},\\3,l'_{3})$, where
$l'_1=a+l_1-6$ and $l'_3=k+l_3-6$. By Lemma \ref{lemma2.2}, we have
\begin{eqnarray*}
b_{2i}(G_{8})&=&b_{2i}(P^{a,k}_{a+k+l_{1}+l_3-3}\cup P^{6}_{7})+b_{2i-2}(P^{a}_{a+l_{1}-2}\cup P^k_{k+l_{3}-2}\cup C_6),\\
b_{2i}(G_{08})&=&b_{2i}(P^{6,6}_{a+k+l_{1}+l_3-3}\cup
P^{6}_{7})+b_{2i-2}(P^{6}_{a+l_{1}-2}\cup P^6_{k+l_{3}-2}\cup C_6).
\end{eqnarray*}
By Proposition \ref{Pro2.3}, we can obtain that $G_8\preceq G_{08}$.

\textbf{Subcase 4.2.} $b=6$, $l_3=2$, $l_1\geq3$ and $l_3\geq3$.

For any values of $l_1$ and $l_3$, let $G_9\in
\Theta_{\uppercase\expandafter{\romannumeral 2}}(n;a,6,k;l_{1},2,l_3)$ and $G_{09}\in
\Theta_{\uppercase\expandafter{\romannumeral 2}}(n;6,6,6;\\l'_{1},2,l'_{3})$, where
$l'_1=a+l_1-6$ and $l'_3=k+l_3-6$. By Lemma \ref{lemma2.2}, we have
\begin{eqnarray*}
b_{2i}(G_{9})&=&b_{2i}(P^{a,k}_{a+k+l_{1}+l_3-3}\cup C_{6})+b_{2i-2}(P^{a}_{a+l_{1}-2}\cup P^{k}_{k+l_3-2}\cup P_{5}),\\
b_{2i}(G_{09})&=&b_{2i}(P^{6,6}_{a+k+l_{1}+l_3-3}\cup
C_{6})+b_{2i-2}(P^{6}_{a+l_{1}-2}\cup P^{6}_{k+l_3-2}\cup P_{5}).
\end{eqnarray*}
By Proposition \ref{Pro2.3}, we have $G_9\preceq G_{09}$.

\textbf{Subcase 4.3.} $b=4$, $l_3=5$, $l_1\geq3$ and $l_3\geq3$.

For any values of $l_1$ and $l_3$, let $G_{10}\in
\Theta_{\uppercase\expandafter{\romannumeral 2}}(n;a,4,k;l_{1},5,l_3)$ and $G_{010}\in
\Theta_{\uppercase\expandafter{\romannumeral 2}}(n;6,6,6;\\ l'_{1},3,l'_{3})$, where
$l'_1=a+l_1-6$ and $l'_3=k+l_3-6$. By Lemma \ref{lemma2.2}, we have
\begin{eqnarray*}
b_{2i}(G_{10})&=&b_{2i}(P^{a,k}_{a+k+l_{1}+l_3-3}\cup P^{4}_{7})+b_{2i-2}(P^{a}_{a+l_{1}-2}\cup P^k_{k+l_{3}-2}\cup P^4_6),\\
b_{2i}(G_{010})&=&b_{2i}(P^{6,6}_{a+k+l_{1}+l_3-3}\cup
P^{6}_{7})+b_{2i-2}(P^{6}_{a+l_{1}-2}\cup P^6_{k+l_{3}-2}\cup C_6).
\end{eqnarray*}
By Proposition \ref{Pro2.3}, we can obtain that $G_{10}\preceq
G_{010}$.

\textbf{Subcase 4.4.} $b=4$, $l_2=4$, $l_1\geq3$ and $l_3\geq3$.

Let $G_{11}\in \Theta_{\uppercase\expandafter{\romannumeral
2}}(n;a,4,k;l_{1},4,l_3)$ and $G_{011}\in
\Theta_{\uppercase\expandafter{\romannumeral
2}}(n;6,6,6;l'_{1},2,l'_{3})$, where $l'_1=a+l_1-6$ and
$l'_3=k+l_3-6$. By Lemma \ref{lemma2.2}, we have
\begin{eqnarray*}
b_{2i}(G_{11})&=&b_{2i}(P^{a,k}_{a+k+l_{1}+l_3-3}\cup P^{4}_{6})+b_{2i-2}(P^{a}_{a+l_{1}-2}\cup P^{k}_{k+l_3-2}\cup P_{5}^{4}),\\
b_{2i}(G_{011})&=&b_{2i}(P^{6,6}_{a+k+l_{1}+l_3-3}\cup
C_{6})+b_{2i-2}(P^{6}_{a+l_{1}-2}\cup P^{6}_{k+l_3-2}\cup P_{5}).
\end{eqnarray*}
Also, $\phi(P^4_5;\lambda)=\lambda^5-3\lambda^3+2\lambda$ and
$\phi(P_5;\lambda)=\lambda^5-4\lambda^3+3\lambda.$ So $P_5^4\prec
P_5$. Then by Proposition \ref{Pro2.3}, we have $G_{11}\preceq
G_{011}$.

\textbf{Case 5.} $\begin{cases}
l_1+a-2\leq 6\\
l_2+b-3\geq 7\\
l_3+k-2\leq 6
\end{cases}$

It is easy to verify that $a\leq6$ and then we have $b\leq k\leq a\leq6$. If $a=b=k=6$, it follows that this lemma holds. Then we focus on other subcases. Without considering graphs with forms $\Gamma_2(\romannumeral 3)$, $\Gamma_2(\romannumeral 4)$, we can distinguish this case into the following three subcases.

\textbf{Subcase 5.1.} $a=k=6$, $b=4$.

It is easy to verify that $l_1=l_3=2$ and $l_2\geq6$. For any
value of $l_2$, let $G_{12}\in
\Theta_{\uppercase\expandafter{\romannumeral 2}}(n;6,4,6;2,l_2,2)$
and $G_{012}\in \Theta_{\uppercase\expandafter{\romannumeral
2}}(n;6,6,6;2,l'_{2},2)$, where $l'_2=l_2-2$. By Lemma
\ref{lemma2.2}, we have
\begin{eqnarray*}
b_{2i}(G_{12})&=&b_{2i}(P^{6,6}_{13}\cup P^{4}_{l_2+2})+b_{2i-2}(C_6 \cup C_6 \cup P^4_{l_{2}+1}),\\
b_{2i}(G_{012})&=&b_{2i}(P^{6,6}_{13}\cup
P^{6}_{l_2+2})+b_{2i-2}(C_6 \cup C_6 \cup P^6_{l_{2}+1}).
\end{eqnarray*}
By Proposition \ref{Pro2.3}, we can obtain that $G_{12}\preceq
G_{012}$.

\textbf{Subcase 5.2.} $a=6$, $k=b=4$, $l_3=4$.

It is easy to verify that $l_1=2$. For fixed $l_2$, let $G_{13}\in
\Theta_{\uppercase\expandafter{\romannumeral 2}}(n;6,4,4;2,l_2,4)$
and $G_{013}\in \Theta_{\uppercase\expandafter{\romannumeral
2}}(n;6,6,6;2,l'_{2},4)$, where $l'_2=l_2-2$. By Lemma
\ref{lemma2.2}, we have
\begin{eqnarray*}
b_{2i}(G_{13})&=&b_{2i}(P^{4,4}_{l_{2}+9}\cup
C_{6})+b_{2i-2}(P_{5}\cup P^{4}_{6}\cup P^{4}_{l_2+2}),\\
b_{2i}(G_{013})&=&b_{2i}(P^{6,6}_{l_{2}+9}\cup
C_{6})+b_{2i-2}(P_{5}\cup C_{6}\cup P^{6}_{l_2+2}).
\end{eqnarray*}
By Proposition \ref{Pro2.3}, we have $G_{13}\preceq G_{013}$.

\textbf{Subcase 5.3.} $a=k=b=4$, $l_1=l_3=4$.

It is easy to verify that $l_2\geq6$. For fixed $l_2$, let
$G_{14}\in \Theta_{\uppercase\expandafter{\romannumeral
2}}(n;4,4,4;4,l_2,4)$ and $G_{014}\in
\Theta_{\uppercase\expandafter{\romannumeral
2}}(n;6,6,6;2,l'_{2},2)$, where $l'_2=l_2-2$. By Lemma
\ref{lemma2.2}, we have
\begin{eqnarray*}
b_{2i}(G_{14})&=&b_{2i}(P^{4,4}_{13}\cup P^{4}_{l_2+2})+b_{2i-2}(P^4_6 \cup P^4_6 \cup P^4_{l_{2}+1}),\\
b_{2i}(G_{014})&=&b_{2i}(P^{6,6}_{13}\cup
P^{6}_{l_2+2})+b_{2i-2}(C_6 \cup C_6 \cup P^6_{l_{2}+1}).
\end{eqnarray*}
Also, $P_6^4 \prec C_6$, and by Proposition \ref{Pro2.3}, we can
obtain that $G_{14}\preceq G_{014}$.

\textbf{Case 6.} $\begin{cases}
l_1+a-2\leq 6\\
l_2+b-3\leq 6\\
l_3+k-2\geq 7
\end{cases}$

It is easy to verify that $a\leq6$ and then we have $b\leq k\leq
a\leq6$. If $a=b=k=6$, it follows that this lemma holds. Then we
focus on other subcases. Without considering graphs with forms
$\Gamma_2(\romannumeral 3)$, $\Gamma_2(\romannumeral 4)$, we can
distinguish this case into the following three subcases.

\textbf{Subcase 6.1.} $a=k=6$, $b=4$, $4 \leq l_2 \leq 5$.

It is easy to verify that $l_1=2$, $l_3 \geq 3$. For any values of $l_2$
and $l_3$, let $G_{15}\in \Theta_{\uppercase\expandafter{\romannumeral
2}}(n;6,4,6;2,l_2,l_{3})$ and $G_{015}\in
\Theta_{\uppercase\expandafter{\romannumeral
2}}(n;6,6,6;2,l'_2,l'_{3})$, where $l'_2=l_2-2$, $l'_3=l_3$.
By Lemma \ref{lemma2.2}, we have
\begin{eqnarray*}
b_{2i}(G_{15})&=&b_{2i}(P^{6,4}_{l_2+9}\cup P^{6}_{l_3+4})+b_{2i-2}(C_6 \cup P^{4}_{l_{2}+2} \cup P^{6}_{l_{3}+3}),\\
b_{2i}(G_{015})&=&b_{2i}(P^{6,6}_{l_2+9}\cup
P^{6}_{l_3+4})+b_{2i-2}(C_6 \cup P^{6}_{l_{2}+2} \cup
P^{6}_{l_{3}+3}).
\end{eqnarray*}
By Proposition \ref{Pro2.3}, we can obtain that $G_{15}\preceq
G_{015}$.

\textbf{Subcase 6.2.} $a=6$, $k=b=4$, $4\leq l_2\leq5$.

It is easy to verify that $l_1=2$, $l_3\geq 5$. For any values of
$l_2$ and $l_3$, let $G_{16}\in
\Theta_{\uppercase\expandafter{\romannumeral
2}}(n;6,4,4;2,l_2,l_{3})$ and $G_{016}\in
\Theta_{\uppercase\expandafter{\romannumeral
2}}(n;6,6,6;2,l'_2,l'_{3})$, where $l'_2=l_2-2$, $l'_3=l_3-2$. By
Lemma \ref{lemma2.2}, we have
\begin{eqnarray*}
b_{2i}(G_{16})&=&b_{2i}(P^{6,4}_{l_{2}+9}\cup P^4_{l_3+2})+b_{2i-2}(C_{6}\cup P^{4}_{l_2+2}\cup P^{4}_{l_3+1}),\\
b_{2i}(G_{016})&=&b_{2i}(P^{6,6}_{l_{2}+9}\cup
P^6_{l_3+2})+b_{2i-2}(C_{6}\cup P^{6}_{l_2+2}\cup P^{6}_{l_3+1}).
\end{eqnarray*}
By Proposition \ref{Pro2.3}, we can obtain that $G_{16}\preceq
G_{016}$.

\textbf{Subcase 6.3.} $a=k=b=4$, $l_1=4$, $4\leq l_2\leq 5$.

It is easy to verify that $l_3\geq 5$. For any values of $l_2$ and
$l_3$, let $G_{17}\in \Theta_{\uppercase\expandafter{\romannumeral
2}}(n;4,4,4;4,l_2,l_{3})$ and $G_{017}\in
\Theta_{\uppercase\expandafter{\romannumeral
2}}(n;6,6,6;2,l'_2,l'_{3})$, where $l'_2=l_2-2$, $l'_3=l_3-2$. By
Lemma \ref{lemma2.2}, we have
\begin{eqnarray*}
b_{2i}(G_{17})&=&b_{2i}(P^{4,4}_{l_2+9}\cup P^{4}_{l_3+2})+b_{2i-2}(P^4_6 \cup P^4_{l_2+2} \cup P^4_{l_{3}+1}),\\
b_{2i}(G_{017})&=&b_{2i}(P^{6,6}_{l_2+9}\cup
P^{6}_{l_3+2})+b_{2i-2}(C_6 \cup P^6_{l_2+2})\cup P^6_{l_{3}+1}).
\end{eqnarray*}
Also, $P^4_6 \prec C_6$, and by Proposition \ref{Pro2.3}, we have
$G_{17}\preceq G_{017}$.

\textbf{Case 7.} $\begin{cases}
l_1+a_1-2\geq 7\\
l_2+b-3\leq 6\\
l_3+k-2\leq 6
\end{cases}$

It is easy to verify that $k\leq6$ and $b\leq6$. If $b\leq k\leq
a\leq6$, with similar analysis in Case $6$ we can obtain that this
lemma holds. Then we consider the case of $a>6\geq k\geq b$. Without
considering graphs with forms $\Gamma_2(\romannumeral 1)$,
$\Gamma_2(\romannumeral 3)$ and $\Gamma_2(\romannumeral 4)$, we can
distinguish this case into the following three subcases.

\textbf{Subcase 7.1.} $k=b=6$, $l_1\geq 3$.

It is easy to verify that $l_3=2$ and $2\leq l_2\leq3$. For any
values of $l_1$ and $l_2$, let $G_{18}\in
\Theta_{\uppercase\expandafter{\romannumeral
2}}(n;a,6,6;l_1,l_2,2)$ and $G_{018}\in
\Theta_{\uppercase\expandafter{\romannumeral
2}}(n;6,6,6;l'_1,l_2,2)$, where $l'_1=a+l_1-6$. By Lemma
\ref{lemma2.2}, we have
\begin{eqnarray*}
b_{2i}(G_{18})&=&b_{2i}(P^{a,6}_{a+l_1+l_{2}+3}\cup C_{6})+b_{2i-2}(P^{a}_{a+l_1-2}\cup P^{6}_{l_2+4}\cup P_{5}),\\
b_{2i}(G_{013})&=&b_{2i}(P^{6,6}_{a+l_1+l_{2}+3}\cup
C_{6})+b_{2i-2}(P^{6}_{a+l_1-2}\cup P^{6}_{l_2+4}\cup P_{5}).
\end{eqnarray*}
By Proposition \ref{Pro2.3}, it follows that $G_{18}\preceq
G_{018}$.

\textbf{Subcase 7.2.} $k=6$, $b=4$, $l_1\geq 3$.

It is easy to verify that $l_3=2$ and $l_2\leq5$. We have $4\leq
l_2\leq5$ since we do not consider graphs with forms
$\Gamma_2(\romannumeral 3)$ and $\Gamma_2(\romannumeral 4)$. For any
values of $l_1$ and $l_3$, let $G_{19}\in
\Theta_{\uppercase\expandafter{\romannumeral
2}}(n;a,4,6;l_1,l_2,2)$ and $G_{019}\in
\Theta_{\uppercase\expandafter{\romannumeral
2}}(n;6,6,6;l'_1,l'_{2},2)$, where $l'_1=a+l_1-6$ and
$l'_2=l_2-2$. By Lemma \ref{lemma2.2}, we have
\begin{eqnarray*}
b_{2i}(G_{19})&=&b_{2i}(P^{a,4}_{a+l_1+l_{2}+1}\cup C_{6})+b_{2i-2}(P^{a}_{a+l_1-2}\cup P^{4}_{l_2+2}\cup P_{5}),\\
b_{2i}(G_{019})&=&b_{2i}(P^{6,6}_{a+l_1+l_{2}+1}\cup
C_{6})+b_{2i-2}(P^{6}_{a+l_1-2}\cup P^{6}_{l_2+2}\cup P_{5}).
\end{eqnarray*}
By Proposition \ref{Pro2.3}, we have $G_{19}\preceq G_{019}$.

\textbf{Subcase 7.3.} $k=b=4$, $l_1\geq 3$, $l_3=4$.

Similar to Subcase 7.2, we have $4\leq l_2\leq5$. Let $G_{20}\in
\Theta_{\uppercase\expandafter{\romannumeral 2}}(n;a,4,4;l_1,l_2,4)$
and $G_{020}\in \Theta_{\uppercase\expandafter{\romannumeral
2}}(n;6,6,6;l'_1,l'_{2},2)$, where $l'_1=a+l_1-6$ and $l'_2=l_2-2$.
By Lemma \ref{lemma2.2}, we have
\begin{eqnarray*}
b_{2i}(G_{20})&=&b_{2i}(P^{a,4}_{a+l_1+l_{2}+1}\cup P^4_{6})+b_{2i-2}(P^{a}_{a+l_1-2}\cup P^{4}_{l_2+2}\cup P_{5}),\\
b_{2i}(G_{020})&=&b_{2i}(P^{6,6}_{a+l_1+l_{2}+1}\cup
C_{6})+b_{2i-2}(P^{6}_{a+l_1-2}\cup P^{6}_{l_2+2}\cup P_{5}).
\end{eqnarray*}
Since $P^4_{6}\prec C_{6}$ and $P^4_{5}\prec P_{5}$, by Proposition
\ref{Pro2.3}, we have $G_{20}\preceq G_{020}$.

\textbf{Case 8.} $\begin{cases}
l_1+a_1-2\leq 6\\
l_2+b-3\leq 6\\
l_3+k-2\leq 6
\end{cases}$

It is easy to verify that $a\leq6$ and then we have $b\leq k\leq
a\leq6$. If $a=b=k=6$, it follows that this lemma holds. Then we
focus on other subcases. Without considering the graphs with forms
$\Gamma_2(\romannumeral 3)$ and $\Gamma_2(\romannumeral 4)$, we can
distinguish this case into the following three subcases.

\textbf{Subcase 8.1.} $a=k=6$, $b=4$.

It is easy to verify that $l_1=l_3=2$. Since $n\geq 20$, we have
$l_2=5$. Let $G_{21}\in \Theta_{\uppercase\expandafter{\romannumeral
2}}(20;6,4,6;2,5,2)$ and $G_{021}\in
\Theta_{\uppercase\expandafter{\romannumeral 2}}(20;6,6,6;2,3,2)$.
By Lemma \ref{lemma2.2}, we have
\begin{eqnarray*}
b_{2i}(G_{21})&=&b_{2i}(P^{6,4}_{14}\cup C_{6})+b_{2i-2}(C_6 \cup P^{4}_{7} \cup P_{5}),\\
b_{2i}(G_{021})&=&b_{2i}(P^{6,6}_{14}\cup C_{6})+b_{2i-2}(C_6 \cup
P^{6}_{7} \cup P_{5}).
\end{eqnarray*}
By Proposition \ref{Pro2.3}, it follows that $G_{21}\preceq
G_{021}$.

\textbf{Subcase 8.2.} $a=6$, $k=b=4$, $l_3=4$.

It is easy to verify that $l_1=2$, $l_2\leq5$. Since $n\geq20$, we
have $l_2=5$. Let $G_{22}\in
\Theta_{\uppercase\expandafter{\romannumeral 2}}(20;6,4,4;2,5,4)$
and $G_{022}\in \Theta_{\uppercase\expandafter{\romannumeral
2}}(20;6,6,6;2,3,2)$. By Lemma \ref{lemma2.2}, we have
\begin{eqnarray*}
b_{2i}(G_{22})&=&b_{2i}(C_{6}\cup P^{4,4}_{14})+b_{2i-2}(P_5\cup P^{4}_{6}\cup P^{4}_{7}),\\
b_{2i}(G_{022})&=&b_{2i}(C_{6}\cup P^{6,6}_{14})+b_{2i-2}(P_5\cup
C_6\cup P^{6}_{7}).
\end{eqnarray*}
By Proposition \ref{Pro2.3}, we have $G_{22}\preceq G_{022}$.

\textbf{Subcase 8.3.} $a=k=b=4$, $l_1=l_3=4$.

It is easy to verify that $l_2\leq5$. Since $n\geq20$, we have
$l_2=5$. Let $G_{23}\in \Theta_{\uppercase\expandafter{\romannumeral
2}}(20;4,4,4;4,5,4)$ and $G_{023}\in
\Theta_{\uppercase\expandafter{\romannumeral 2}}(20;6,6,6;2,3,2)$.
By Lemma \ref{lemma2.2}, we have
\begin{eqnarray*}
b_{2i}(G_{23})&=&b_{2i}(P^{4,4}_{13}\cup P^{4}_{7})+b_{2i-2}(P^4_6 \cup P^4_{6} \cup P^4_{6}),\\
b_{2i}(G_{023})&=&b_{2i}(P^{6,6}_{13}\cup P^{6}_{7})+b_{2i-2}(C_6
\cup C_6\cup C_6).
\end{eqnarray*}
Since $P^4_6 \prec C_6$, and by Proposition \ref{Pro2.3}, we can
obtain that $G_{23}\preceq G_{023}$.

The proof is now complete. \qed

\begin{figure}[h,t,b,p]
\begin{center}
\includegraphics[scale = 0.6]{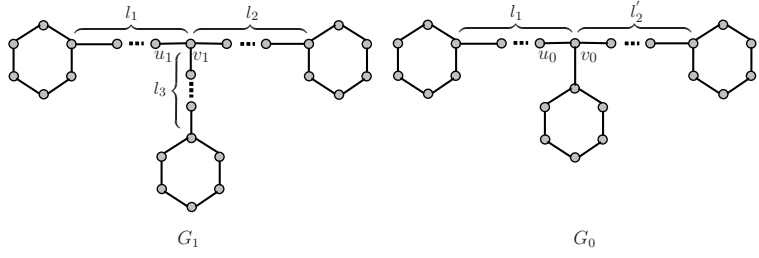}
\caption{Graphs for Lemma \ref{lemma3.5}.} \label{fig8}
\end{center}
\end{figure}

\begin{lem}\label{lemma3.5}
For any graph $G\in \Theta_{\uppercase\expandafter{\romannumeral
2}}(n;6,6,6;l_{1},l_{2},l_{3})$, there exists a graph $H\in
\Theta_{\uppercase\expandafter{\romannumeral
2}}(n;6,6,6;l'_{1},l'_{2},2)$ such that $G\prec H$.
\end{lem}
\pf For fixed parameters $n$, $l_{1}$, $l_{2}$ and $l_{3}$, let
$G_{1}\in\Theta_{\uppercase\expandafter{\romannumeral
2}}(n;6,6,6;l_{1},l_{2},l_{3})$ and
$G_{0}\in\Theta_{\uppercase\expandafter{\romannumeral
2}}(n;6,6,6;l_{1},l'_{2},2)$ (as shown in Figure \ref{fig8}). It is
easy to verify that $l'_{2}=l_{2}+l_{3}-2$ and it suffices to show
that $G_{1}\prec G_{0}$.

By Lemma \ref{lemma2.2} we have
\begin{eqnarray*}
b_{2i}(G_{1})&=&b_{2i}(G_{1}-u_{1}v_{1})+b_{2i-2}(G_{1}-u_{1}-v_{1})\\
&=&b_{2i}(G_{1}-u_{1}v_{1})+b_{2i-2}(P_{l_{1}+3}^{6}\cup P_{l_{2}+4}^{6}\cup P^{6}_{l_{3}+4})\\
&=&b_{2i}(G_{1}-u_{1}v_{1})+b_{2i-2}(P_{l_{1}+3}^{6}\cup P_{l_{2}+4}^{6}\cup C_{6}\cup P_{l_{3}-2})\\
&&+b_{2i-4}(P_{l_{1}+3}^{6}\cup P_{l_{2}+4}^{6}\cup P_{l_{3}-3}\cup
P_{5}),\\
b_{2i}(G_{0})&=&b_{2i}(G_{0}-u_{0}v_{0})+b_{2i-2}(G_{0}-u_{0}-v_{0})\\
&=&b_{2i}(G_{0}-u_{0}v_{0})+b_{2i-2}(P_{l_{1}+3}^{6}\cup P_{l_{2}+l_{3}+2}^{6}\cup C_{6})\\
&=&b_{2i}(G_{0}-u_{0}v_{0})+b_{2i-2}(P_{l_{1}+3}^{6}\cup P_{l_{2}+4}^{6}\cup C_{6}\cup P_{l_{3}-2})\\
&&+b_{2i-4}(P_{l_{1}+3}^{6}\cup P_{l_{2}+3}^{6}\cup P_{l_{3}-3}\cup
C_{6}).
\end{eqnarray*}
Since $b_{2i}(G_{1}-u_{1}v_{1})=b_{2i}(G_{0}-u_{0}v_{0})$, then we
only need to compare $b_{2j}(P_{l_{2}+4}^{6}\cup P_{5})$ and
$b_{2j}(P_{l_{2}+3}^{6}\cup C_{6})$. Also by Lemma \ref{lemma2.2} we
have
\begin{eqnarray*}
b_{2j}(P_{l_{2}+4}^{6}\cup P_{5})&=&b_{2j}(P_{l_{2}+3}^{6}\cup P_{5}\cup P_{1})+b_{2j-2}(P_{l_{2}+2}^{6}\cup P_{5})\\
&=&b_{2j}(P_{l_{2}+3}^{6}\cup P_{5}\cup P_{1})+b_{2j-2}(P_{l_{2}+2}^{6}\cup P_{4}\cup P_{1})\\
&&+b_{2j-4}(P_{l_{2}+2}^{6}\cup P_{3})\\
&=&b_{2j}(P_{l_{2}+3}^{6}\cup P_{5}\cup P_{1})+b_{2j-2}(P_{l_{2}+2}^{6}\cup P_{4}\cup P_{1})\\
&&+b_{2j-4}(P_{l_{2}+1}^{6}\cup P_{3}\cup P_{1})
+b_{2j-6}(P_{l_{2}}^{6}\cup P_{3})\\
&=&b_{2j}(P_{l_{2}+3}^{6}\cup P_{5}\cup P_{1})+b_{2j-2}(P_{l_{2}+2}^{6}\cup P_{4}\cup P_{1})\\
&&+b_{2j-4}(P_{l_{2}+1}^{6}\cup P_{3}\cup P_{1})
+b_{2j-6}(C_{6}\cup P_{l_{2}-6}\cup P_{3})\\
&&+b_{2j-8}(P_{5}\cup P_{l_{2}-7}\cup P_{3}),
\end{eqnarray*}
and \begin{eqnarray*}
b_{2j}(P_{l_{2}+3}^{6}\cup C_{6})&=&b_{2j}(P_{l_{2}+3}^{6}\cup P_{6})+b_{2j-2}(P_{l_{2}+3}^{6}\cup P_{4})+2b_{2j-6}(P_{l_{2}+3}^{6})\\
&=&b_{2j}(P_{l_{2}+3}^{6}\cup P_{6})+b_{2j-2}(P_{l_{2}+2}^{6}\cup P_{4}\cup P_{1})\\
&&+b_{2j-4}(P_{l_{2}+1}^{6}\cup P_{4})+2b_{2j-6}(P_{l_{2}+3}^{6})\\
&=&b_{2j}(P_{l_{2}+3}^{6}\cup P_{6})+b_{2j-2}(P_{l_{2}+2}^{6}\cup P_{4}\cup P_{1})\\
&&+b_{2j-4}(P_{l_{2}+1}^{6}\cup P_{3}\cup P_{1})+b_{2j-6}(P_{l_{2}+1}^{6}\cup P_{2})+2b_{2j-6}(P_{l_{2}+3}^{6})\\
&=&b_{2j}(P_{l_{2}+3}^{6}\cup P_{6})+b_{2j-2}(P_{l_{2}+2}^{6}\cup P_{4}\cup P_{1})\\
&&+b_{2j-4}(P_{l_{2}+1}^{6}\cup P_{3}\cup P_{1})+b_{2j-6}(C_{6}\cup P_{l_{2}-5}\cup P_{2})\\
&&+b_{2j-8}(P_{5}\cup P_{l_{2}-6}\cup
P_{2})+2b_{2j-6}(P_{l_{2}+3}^{6}).
\end{eqnarray*}
By Lemma \ref{lemma2.4} and Proposition \ref{Pro2.3} we have
$P_{l_{2}+4}^{6}\cup P_{5}\prec P_{l_{2}+3}^{6}\cup C_{6}$. Also
consider Proposition \ref{Pro2.3}, we can obtain that $G_{1}\prec
G_{0}$. \qed

\begin{figure}[h,t,b,p]
\begin{center}
\includegraphics[scale = 0.7]{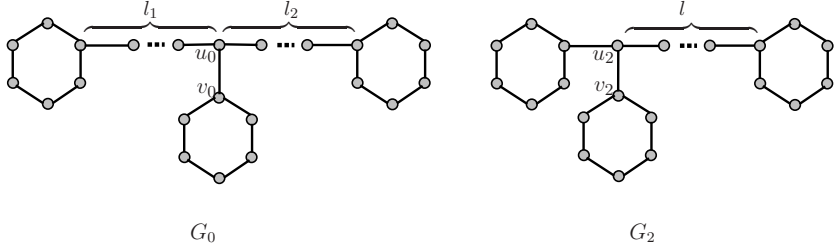}
\caption{Graphs for Lemma 3.6} \label{fig9}
\end{center}
\end{figure}

\begin{lem}\label{lemma3.6}
For any graph $G\in \Theta_{\uppercase\expandafter{\romannumeral 2}}(n;6,6,6;l_{1},l_{2},2)$, there exists a graph $H\in \Theta_{\uppercase\expandafter{\romannumeral 2}}(n;6,6,6;l,2,2)$ such that $G\prec H$.
\end{lem}
\pf For fixed parameters $n$, $l_{1}$ and $l_{2}$, let
$G_{0}\in\Theta_{\uppercase\expandafter{\romannumeral
2}}(n;6,6,6;l_{1},l_{2},2)$ and
$G_{2}\in\Theta_{\uppercase\expandafter{\romannumeral
2}}(n;6,6,6;l,2,2)$ (as shown in Figure \ref{fig9}). It is easy to
verify that $l=l_{1}+l_{2}-2$ and it suffices to show that
$G_{0}\prec G_{2}$.

By Lemma \ref{lemma2.2} we have
\begin{eqnarray*}
b_{2i}(G_{0})&=&b_{2i}(G_{0}-u_{0}v_{0})+b_{2i-2}(G_{0}-u_{0}-v_{0})\\ &=&b_{2i}(G_{0}-u_{0}v_{0})+b_{2i-2}(P_{l_{1}+4}^{6}\cup P_{l_{2}+4}^{6}\cup P_{5})\\
&=&b_{2i}(G_{0}-u_{0}v_{0})+b_{2i-2}(P_{l_{2}+4}^{6}\cup C_{6}\cup P_{l_{1}-2}\cup P_{5})\\
&&+b_{2i-4}(P_{l_{2}+4}^{6}\cup P_{5}\cup P_{l_{1}-3}\cup
P_{5}),\\
b_{2i}(G_{2})&=&b_{2i}(G_{2}-u_{2}v_{2})+b_{2i-2}(G_{2}-u_{2}-v_{2})\\
&=&b_{2i}(G_{2}-u_{2}v_{2})+b_{2i-2}(P_{l_{1}+l_{2}+2}^{6}\cup C_{6}\cup P_{5})\\
&=&b_{2i}(G_{2}-u_{2}v_{2})+b_{2i-2}(P_{l_{2}+4}^{6}\cup C_{6}\cup P_{l_{1}-2}\cup P_{5})\\
&&+b_{2i-4}(P_{l_{2}+3}^{6}\cup C_{6}\cup P_{l_{1}-3}\cup P_{5}).
\end{eqnarray*}
Since $b_{2i}(G_{0}-u_{0}v_{0})=b_{2i}(G_{2}-u_{2}v_{2})$, then we
only need to compare
$b_{2j}(P_{l_{2}+4}^{6}\cup P_{5})$ with $b_{2j}(P_{l_{2}+3}^{6}\cup
C_{6})$. With similar analysis in Lemma \ref{lemma3.5}, we can
obtain that $G_{0}\prec G_{2}$. \qed

From Theorem \ref{Thm3.4}, Lemmas \ref{lemma3.5} and \ref{lemma3.6},
we can easily obtain the following result.

\begin{thm}\label{Thm3.7}
For any graph $G\in \Theta_{\uppercase\expandafter{\romannumeral
2}}(n;a,b,k;l_{1},l_{2},l_{3})$, if $G$ is not an element of the
special graph class $\Gamma_{2}$, then there exists a graph $H\in
\Theta_{\uppercase\expandafter{\romannumeral 2}}(n;6,6,6;n-17,2,2)$
such that $G\preceq H$, and the equality holds if and only if
$G\cong H$.
\end{thm}

\begin{thm}\label{Thm3.8}
For any graph $G\in \Theta_{\uppercase\expandafter{\romannumeral
1}}(n;a,b,k;l_{1},l_{2};2)\setminus \Gamma_{1}$, there exists a
graph $H\in \Theta_{\uppercase\expandafter{\romannumeral
1}}(n;6,6,6;l'_{1},l'_{2};2)$ such that $G\preceq H$.
\end{thm}
\pf Without loss of generality, we may assume that $l_1\geq l_2$. We
will discuss the following four cases.

\textbf{Case 1.} $\begin{cases}
l_1+a-1\geq 9\\
l_2+k-1\geq 8
\end{cases}$

Considering the values of $l_1$ and $l_2$, we distinguish this case into the following four subcases.

\textbf{Subcase 1.1.} $l_1\geq4$.

For any values of $l_1$ and $l_2$, let $G_1\in
\Theta_{\uppercase\expandafter{\romannumeral
1}}(n;a,b,k;l_{1},l_{2};2)$ and $G_{01}\in
\Theta_{\uppercase\expandafter{\romannumeral
1}}(n;6,6,6;l'_{1},l'_{2};2)$, where $l'_{1}=a+l_1-6$. By
lemma \ref{lemma2.2}, we have
\begin{eqnarray*}
b_{2i}(G_{1})&=&b_{2i}(P^{a}_{a+l_{1}-2}\cup P^{b,k}_{b+k+l_2-2})+b_{2i-2}(P^{a}_{a+l_{1}-3}\cup P^{k}_{b+k+l_2-3}),\\
b_{2i}(G_{01})&=&b_{2i}(P^{6}_{a+l_{1}-2}\cup
P^{6,6}_{b+k+l_2-2})+b_{2i-2}(P^{6}_{a+l_{1}-3}\cup
P^{6}_{b+k+l_2-3}).
\end{eqnarray*}
By Proposition \ref{Pro2.3}, we can obtain that $G_1\preceq G_{01}$.

\textbf{Subcase 1.2.} $l_1=l_2=3$ and $b\geq6$.

It is easy to verify that $a\geq8$ and $k\geq6$. Let $G_2\in
\Theta_{\uppercase\expandafter{\romannumeral 1}}(n;a,b,k;3,3;2)$ and
$G_{02}\in \Theta_{\uppercase\expandafter{\romannumeral
1}}(n;6,6,6;l'_{1},l'_{2};2)$, where $l'_{1}=a-3$ and
$l'_{2}=b+k-9$. By lemma \ref{lemma2.2}, we have
\begin{eqnarray*}
b_{2i}(G_{2})&=&b_{2i}(P^{a,k}_{n})+b_{2i-2}(P^{a}_{a+1}\cup P^{k}_{k+1}\cup P_{b-2})+(-1)^{1+\frac{b}{2}}2b_{2i-b}(P^{a}_{a+1}\cup P^{k}_{k+1}),\\
b_{2i}(G_{02})&=&b_{2i}(P^{6,6}_{n})+b_{2i-2}(P^{6}_{a+1}\cup
P^{6}_{b+k-5}\cup P_{4})+2b_{2i-6}(P^{6}_{a+1}\cup P^{6}_{b+k-5}).
\end{eqnarray*}
Then we compare
$b_{2j}(P^{k}_{k+1}\cup P_{b-2})$ with $b_{2j}(P^{6}_{b+k-5}\cup
P_{4})$. By Lemma \ref{lemma2.2} we have
\begin{eqnarray*}
b_{2j}(P^{k}_{k+1}\cup P_{b-2})&=&b_{2j}(P_{k+1}\cup P_{b-2})+b_{2j-2}(P_{k-2}\cup P_{b-2}\cup P_{1}),\\
b_{2j}(P^{6}_{b+k-5}\cup P_{4})&=&b_{2j}(P_{b+k-5}\cup
P_{4})+b_{2j-2}(P_{b+k-11}\cup P_{4}\cup P_{4}).
\end{eqnarray*}
Since $b\geq6$ and $k\geq6$, by Lemma \ref{lemma2.4}, we have
$P_{k+1}\cup P_{b-2}\prec P_{b+k-5}\cup P_{4}$ and $P_{k-2}\cup
P_{b-2}\cup P_{1}\prec P_{k-2}\cup P_{b-5}\cup P_{4}\preceq
P_{b+k-11}\cup P_{4}\cup P_{4}$. Then we can obtain that $
P_{k+1}^{k}\cup P_{b-2}\preceq P_{b+k-5}^{6}\cup P_{4}$. Also, since
$b\geq6$, then $k+1\leq b+k-5$, by Proposition \ref{Pro2.3} we have
$P_{k+1}^{k}\preceq P_{k+1}^{6}\preceq P_{b+k-5}^{6}$. Also by
Proposition \ref{Pro2.3}, we can obtain that $G_2\preceq G_{02}$.

\textbf{Subcase 1.3.} $l_1=l_2=3$, $b=4$ and $k=6$.

It is easy to verify that $a\geq8$. Let $G_3\in
\Theta_{\uppercase\expandafter{\romannumeral 1}}(n;a,4,6;3,3;2)$ and
$G_{03}\in\Theta_{\uppercase\expandafter{\romannumeral
1}}(n;6,6,6;l'_{1},3;2)$, where $l'_{1}=a-5$. By Lemma
\ref{lemma2.2}, we have
\begin{eqnarray*}
b_{2i}(G_{3})&=&b_{2i}(P^{a,4}_{a+5}\cup P_{7}^{6})+b_{2i-2}(P^{a}_{a+4}\cup C_6),\\
b_{2i}(G_{03})&=&b_{2i}(P^{6,6}_{a+5}\cup
P_{7}^{6})+b_{2i-2}(P^{6}_{a+4}\cup C_6).
\end{eqnarray*}
From Proposition \ref{Pro2.3}, it follows that $G_3\preceq G_{03}$.

\textbf{Subcase 1.4.} $l_1=l_2=3$, $b=4$, $k\geq8$ or $l_1=3$,
$l_2=2$ or $l_1=l_2=2$.

The graphs in this case belong to $\Gamma_1(\romannumeral 1)$ or
$\Gamma_1(\romannumeral 2)$, so we do not consider them.

\textbf{Case 2.} $\begin{cases}
l_1+a-1\leq 8\\
l_2+k-1\geq 8
\end{cases}$

It is easy to verify that $a\leq6$. Without considering graphs of
form $\Gamma_1(\romannumeral 3)$, we distinguish this case into the
following two subcases.

\textbf{Subcase 2.1.} $a=6$.

It is easy to verify that $l_1=2\ or\ 3$. If $l_1=3$, $l_2=3$, then
let $G_4\in \Theta_{\uppercase\expandafter{\romannumeral
1}}(n;6,b,k;3,3;2)$ and $G_{04}=
\Theta_{\uppercase\expandafter{\romannumeral
1}}(n;6,6,6;3,l'_{2};2)$, where $l'_{2}=b+k-9$. By Lemma
\ref{lemma2.2}, we have
\begin{eqnarray*}
b_{2i}(G_{4})&=&b_{2i}(P^{b,k}_{b+k+1}\cup P^{6}_{7} )+b_{2i-2}(P^{k}_{b+k}\cup C_6),\\
b_{2i}(G_{04})&=&b_{2i}(P^{6,6}_{b+k+1}\cup P^{6}_{7}
)+b_{2i-2}(P^{6}_{b+k}\cup C_6).
\end{eqnarray*}
By Proposition \ref{Pro2.3}, we have $G_4\preceq G_{04}$.

If $l_1=3$, $l_2=2$, then let $G_5\in
\Theta_{\uppercase\expandafter{\romannumeral 1}}(n;6,b,k;3,2;2)$ and
$G_{05}= \Theta_{\uppercase\expandafter{\romannumeral
1}}(n;6,6,6;3,l''_{2};2)$, where $l''_{2}=b+k-10$. With
similar analysis, it follows that $G_5\preceq G_{05}$. If
$l_1=l_2=2$, then let $G_6\in
\Theta_{\uppercase\expandafter{\romannumeral 1}}(n;6,b,k;2,2;2)$ and
$G_{06}\in \Theta_{\uppercase\expandafter{\romannumeral
1}}(n;6,6,6;2,l'''_{2};2)$, where $l'''_{2}=b+k-10$. With
similar analysis, we can obtain that $G_6\preceq G_{06}$.

\textbf{Subcase 2.2.} $a=4$:

It is easy to verify that $l_1\leq5$. Since we do not consider
graphs with form $\Gamma_1(\romannumeral 3)$, we have $4\leq
l_1\leq5$. If $l_1=5$, let $G_7\in
\Theta_{\uppercase\expandafter{\romannumeral 1}}(n;4,b,k;5,l_2;2)$
and $G_{07}\in \Theta_{\uppercase\expandafter{\romannumeral
1}}(n;6,6,6;3,l'_{2};2)$, where $l'_{2}=b+k+l_2-12$. By Lemma
\ref{lemma2.2}, we have
\begin{eqnarray*}
b_{2i}(G_{7})&=&b_{2i}(P^{b,k}_{b+k+l_2-2}\cup P^{4}_{7} )+b_{2i-2}(P^{k}_{b+k+l_2-3}\cup P^{4}_{6}),\\
b_{2i}(G_{07})&=&b_{2i}(P^{6,6}_{b+k+l_2-2}\cup P^{6}_{7}
)+b_{2i-2}(P^{6}_{b+k+l_2-3}\cup C_{6}).
\end{eqnarray*}
From Proposition \ref{Pro2.3}, it follows that $G_7\preceq G_{07}$.
If $l_1=4$, let $G_8\in \Theta_{\uppercase\expandafter{\romannumeral
1}}(n;4,b,k;4,l_2;2)$ and $G_{08}\in
\Theta_{\uppercase\expandafter{\romannumeral
1}}(n;6,6,6;2,l'_{2};2)$, where $l'_{2}=b+k+l_2-12$. By Lemma
\ref{lemma2.2}, we have
\begin{eqnarray*}
b_{2i}(G_{8})&=&b_{2i}(P^{b,k}_{b+k+l_2-2}\cup P^{4}_{6} )+b_{2i-2}(P^{k}_{b+k+l_2-3}\cup P^{4}_{5}),\\
b_{2i}(G_{08})&=&b_{2i}(P^{6,6}_{b+k+l_2-2}\cup C_{6}
)+b_{2i-2}(P^{6}_{b+k+l_2-3}\cup P_{5}).
\end{eqnarray*}
Since $P^{4}_{5}\prec P_5$, then from Proposition \ref{Pro2.3}, it
follows that $G_8\preceq G_{08}$.

\textbf{Case 3.} $\begin{cases}
l_1+a-1\geq 9\\
l_2+k-1\leq 7
\end{cases}$

Without considering graphs with form $\Gamma_1(\romannumeral 4)$, we
distinguish this case into the following two subcases.

\textbf{Subcase 3.1.} $k=6$.

It is easy to verify that $l_2=2$. For any value of $l_1$, let
$G_9\in \Theta_{\uppercase\expandafter{\romannumeral
1}}(n;a,b,6;l_1,2;2)$ and $G_{09}\in
\Theta_{\uppercase\expandafter{\romannumeral
1}}(n;6,6,6;l'_{1},2;2)$, where $l'_{1}=a+b+l_1-12$. By Lemma
\ref{lemma2.2}, we have
\begin{eqnarray*}
b_{2i}(G_{9})&=&b_{2i}(P^{a,b}_{a+b+l_1-2}\cup C_6)+b_{2i-2}(P^{a}_{a+b+l_1-3}\cup P_{5}),\\
b_{2i}(G_{09})&=&b_{2i}(P^{6,6}_{a+b+l_1-2}\cup
C_6)+b_{2i-2}(P^{6}_{a+b+l_1-3}\cup P_{5}).
\end{eqnarray*}
By Proposition \ref{Pro2.3}, we can obtain that $G_9\preceq G_{09}$.

\textbf{Subcase 3.2.} $k=4$.

It is easy to verify that $l_2\leq4$. Since we do not consider
graphs with form $\Gamma_1(\romannumeral 4)$, we have $l_2=4$. For
any value of $l_1$, let $G_{10}\in
\Theta_{\uppercase\expandafter{\romannumeral 1}}(n;a,b,4;l_1,4;2)$
and $G_{010}\in \Theta_{\uppercase\expandafter{\romannumeral
1}}(n;6,6,6;l'_{1},2;2)$, where $l'_{1}=a+b+l_1-12$. By Lemma
\ref{lemma2.2}, we have
\begin{eqnarray*}
b_{2i}(G_{10})&=&b_{2i}(P^{a,b}_{a+b+l_1-2}\cup P^4_6)+b_{2i-2}(P^{a}_{a+b+l_1-3}\cup P^4_{5}),\\
b_{2i}(G_{010})&=&b_{2i}(P^{6,6}_{a+b+l_1-2}\cup
C_6)+b_{2i-2}(P^{6}_{a+b+l_1-3}\cup P_{5}).
\end{eqnarray*}
Since $P^{4}_{5}\prec P_5$, by Proposition \ref{Pro2.3}, we can
obtain that $G_{10}\preceq G_{010}$.

\textbf{Case 4.} $\begin{cases}
l_1+a-1\leq 8\\
l_2+k-1\leq 7
\end{cases}$

It is easy to verify that $a\leq6$ and $k\leq6$. Without considering
graphs with form $\Gamma_1(\romannumeral 5)$, we distinguish this
case into the following two subcases.

\textbf{Subcase 4.1.} $a=6$:

It is easy to verify that $l_1\leq3$. If $l_1=3$, then let
$G_{11}\in\Theta_{\uppercase\expandafter{\romannumeral
1}}(n;6,b,k;3,l_2;2)$ and $G_{011}=
\Theta_{\uppercase\expandafter{\romannumeral
1}}(n;6,6,6;3,l'_{2};2)$, where $l'_{2}=b+k+l_2-12$. By Lemma
\ref{lemma2.2}, we have
\begin{eqnarray*}
b_{2i}(G_{11})&=&b_{2i}(P^{b,k}_{b+k+l_2-2}\cup P^{6}_{7} )+b_{2i-2}(P^{k}_{b+k+l_2-3}\cup C_6),\\
b_{2i}(G_{011})&=&b_{2i}(P^{6,6}_{b+k+l_2-2}\cup P^{6}_{7}
)+b_{2i-2}(P^{6}_{b+k+l_2-3}\cup C_6).
\end{eqnarray*}
By Proposition \ref{Pro2.3}, we can obtain that $G_{11}\preceq
G_{011}$.

If $l_1=2$, since $l_1\geq l_2$, we have $l_2=2$. Let $G_{12}\in
\Theta_{\uppercase\expandafter{\romannumeral 1}}(n;6,b,k;2,2;2)$ and
$G_{012}\in \Theta_{\uppercase\expandafter{\romannumeral
1}}(n;6,6,6;2,l'_{2};2)$, where $l'_{2}=b+k-10$. With similar
analysis, it follows that $G_{12}\preceq G_{012}$.

\textbf{Subcase 4.2.} $a=4$.

It is easy to verify that $l_1\leq5$. Since we do not consider graphs with form
$\Gamma_1(\romannumeral 4)$, then we have $4\leq l_1\leq5$. If $l_1=5$ , then let
$G_{13}\in \Theta_{\uppercase\expandafter{\romannumeral 1}}(n;4,b,k;5,l_2;2)$ and
$G_{013}\in \Theta_{\uppercase\expandafter{\romannumeral 1}}(n;6,6,6;3,l'_{2};2)$, where
$l'_{2}=b+k+l_2-12$. By Lemma \ref{lemma2.2}, we have
$b_{2i}(G_{13})=b_{2i}(P^{b,k}_{b+k+l_2-2}\cup P^{4}_{7}
)+b_{2i-2}(P^{k}_{b+k+l_2-3}\cup P^{4}_{6})$ and
$b_{2i}(G_{013})=b_{2i}(P^{6,6}_{b+k+l_2-2}\cup P^{6}_{7}
)+b_{2i-2}(P^{6}_{b+k+l_2-3}\cup C_{6})$. By Proposition \ref{Pro2.3}, we can obtain
that $G_{13}\preceq G_{013}$.

If $l_1=4$ , then let $G_{14}\in
\Theta_{\uppercase\expandafter{\romannumeral 1}}(n;4,b,k;4,l_2;2)$
and $G_{014}= \Theta_{\uppercase\expandafter{\romannumeral
1}}(n;6,6,6;2,l'_{2};2)$, where $l'_{2}=b+k+l_2-12$. With
similar analysis we can obtain that $G_{14}\preceq G_{014}$.

The proof is thus complete.\qed

\begin{figure}[h,t,b,p]
\begin{center}
\includegraphics[scale = 0.6]{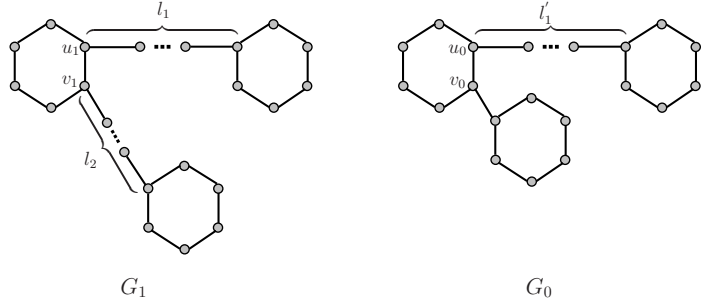}
\caption{Graphs for Lemma \ref{lemma3.9}.} \label{fig10}
\end{center}
\end{figure}

\begin{lem}\label{lemma3.9}
For any graph $G\in \Theta_{\uppercase\expandafter{\romannumeral
1}}(n;6,6,6;l_{1},l_{2};2)$, there exists a graph $H\in
\Theta_{\uppercase\expandafter{\romannumeral
1}}(n;6,6,\\ 6;l'_{1},2;2)$ such that $G\prec H$.
\end{lem}
\pf For fixed parameters $n$, $a$, $b$, $k$, $l_{1}$ and $l_{2}$,
let $G_{1}\in\Theta_{\uppercase\expandafter{\romannumeral
1}}(n;a,b,k;l_{1},l_{2};2)$ and
$G_{0}=\Theta_{\uppercase\expandafter{\romannumeral
1}}(n;a,b,k;l'_{1},2;2)$ (as shown in Figure \ref{fig10}). It is
easy to verify that $l'_{1}=l_{1}+l_{2}-2$ and it suffices to show
that $G_{1}\prec G_{0}$.

By Lemma \ref{lemma2.2} we have
\begin{eqnarray*}
b_{2i}(G_{1})&=&b_{2i}(G_{1}-u_{1}v_{1})+b_{2i-2}(G_{1}-u_{1}-v_{1})
+2b_{2i-6}(G_{1}-C_{6})\\
&=&b_{2i}(G_{1}-u_{1}v_{1})+b_{2i-2}(P_{l_{1}+4}^{6}\cup
P_{l_{2}+4}^{6}\cup P_{4})+2b_{2i-6}(P_{l_{1}+4}^{6}\cup
P_{l_{2}+4}^{6}),\\
b_{2i}(G_{0})&=&b_{2i}(G_{0}-u_{0}v_{0})+b_{2i-2}(G_{2}-u_{0}-v_{0})
+2b_{2i-6}(G_{0}-C_{6})\\
&=&b_{2i}(G_{0}-u_{0}v_{0})+b_{2i-2}(P_{l'_{1}+4}^{6}\cup
C_{6}\cup P_{4})+2b_{2i-6}(P_{l'_{1}+4}^{6}\cup C_{6}).
\end{eqnarray*}
Since $b_{2i}(G_{1}-u_{1}v_{1})=b_{2i}(G_{0}-u_{0}v_{0})$, and
considering Proposition \ref{Pro2.3}, we try to compare
$b_{2j}(P_{l_{1}+4}^{6}\cup P_{l_{2}+4}^{6})$ with
$b_{2j}(P_{l'_{1}+4}^{6}\cup C_{6})$. Also by Lemma \ref{lemma2.2}
we have
\begin{eqnarray*}
b_{2j}(P_{l_{1}+4}^{6}\cup P_{l_{2}+4}^{6})
&=&b_{2j}(P_{l_{1}+4}^{6}\cup C_{6}\cup P_{l_{2}-2})+b_{2j-2}(P_{l_{1}+4}^{6}
\cup P_{l_{2}-3}\cup P_{5}),\\
b_{2j}(P_{l'_{1}+4}^{6}\cup C_{6}) &=&b_{2j}(P_{l_{1}+4}^{6}\cup
C_{6}\cup P_{l_{2}-2})+b_{2j-2}(P_{l_{1}+3}^{6}\cup P_{l_{2}-3}\cup
C_{6}).
\end{eqnarray*}
With similar analysis in Lemma \ref{lemma3.5}, we have
$P_{l_{1}+4}^{6}\cup P_{5}\prec P_{l_{1}+3}^{6}\cup C_{6}$. Applying
Proposition \ref{Pro2.3}, we can obtain $G_{1}\prec G_{0}$. \qed

\begin{figure}[h,t,b,p]
\begin{center}
\includegraphics[scale = 0.6]{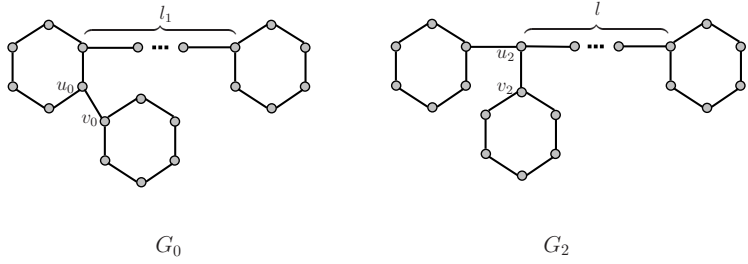}
\caption{Graphs for Lemma \ref{lemma3.10}.} \label{fig11}
\end{center}
\end{figure}

\begin{lem}\label{lemma3.10}
For any graph $G\in \Theta_{\uppercase\expandafter{\romannumeral
1}}(n;6,6,6;l_{1},2;2)$, there exists a graph $H\in
\Theta_{\uppercase\expandafter{\romannumeral 2}}(n;6,6,6;l,2,2)$
such that $\mathcal {E}(G)<\mathcal {E}(H)$.
\end{lem}
\pf For fixed parameters $l_{1}$ and $l$, let
$G_{0}\in\Theta_{\uppercase\expandafter{\romannumeral 1}}(n;6,6,6;l_{1},2;2)$ and
$G_{2}\in\Theta_{\uppercase\expandafter{\romannumeral 2}}(n;6,6,\\6;l,2,2)$ (as shown in
Figure \ref{fig11}), where $l=l_1-1$, i.e., $l_1=l+1$. It suffices to show that
$G_{0}\prec G_{2}$.
 \\
By Lemma \ref{lemma2.2} we have
\begin{eqnarray*}
b_{2i}(G_{0})&=&b_{2i}(G_{0}-u_{0}v_{0})+b_{2i-2}(G_{0}-u_{0}-v_{0})\\
&=&b_{2i}(G_{0}-u_{0}v_{0})+b_{2i-2}(P_{l+10}^{6}\cup P_{5})\\
&=&b_{2i}(G_{0}-u_{0}v_{0})+b_{2i-2}(P_{l+4}^{6}\cup P_{6}\cup P_{5})
+b_{2i-4}(P_{l+3}^{6}\cup P_{5}\cup P_{5}),\\
b_{2i}(G_{2})&=&b_{2i}(G_{2}-u_{2}v_{2})+b_{2i-2}(G_{2}-u_{2}-v_{2})\\
&=&b_{2i}(G_{2}-u_{2}v_{2})+b_{2i-2}(P_{l+4}^{6}\cup C_{6}\cup P_{5})\\
&=&b_{2i}(G_{2}-u_{2}v_{2})+b_{2i-2}(P_{l+4}^{6}\cup P_{6}\cup
P_{5})+b_{2i-4}(P_{l+4}^{6}\cup P_{4}\cup P_{5}).
\end{eqnarray*}
Since $b_{2i}(G_{0}-u_{0}v_{0})=b_{2i}(G_{2}-u_{2}v_{2})$, then we only need to verify $
P_{l+3}^{6}\cup P_{5}\cup P_{5}\prec P_{l+4}^{6}\cup P_{4}\cup P_{5}. $ By Lemma
\ref{lemma2.2}, we have
\begin{eqnarray*}
b_{2i}(P_{l+3}^{6}\cup P_{5})
&=&b_{2i}(P_{l+3}\cup P_{5})+b_{2i-2}(P_{l-3}
\cup P_{5}\cup P_{4})+2b_{2i-6}(P_{l-3}
\cup P_{5}),\\
b_{2i}(P_{l+4}^{6}\cup P_{4}) &=&b_{2i}(P_{l+4}\cup
P_{4})+b_{2i-2}(P_{l-2} \cup P_{4}\cup P_{4})+2b_{2i-6}(P_{l-2} \cup
P_{4}).
\end{eqnarray*}
From Lemma \ref{lemma2.4}, we can obtain that $P_{l+3}\cup P_5\prec
P_{l+4}\cup P_4$ and if $l\neq5$, $P_{l-3}\cup P_5\prec P_{l-2}\cup
P_4$, then $P_{l-3}\cup P_5\cup P_4\prec P_{l-2}\cup P_4\cup P_4$.
So from Proposition \ref{Pro2.3}, it follows that $P_{l+3}^6\cup
P_5\prec P_{l+4}^6\cup P_4$ and then $G_0\prec G_2$. If $l=5$, then
$G_{0}\in\Theta_{\uppercase\expandafter{\romannumeral
1}}(22;6,6,6;6,2;2)$ and
$G_{2}\in\Theta_{\uppercase\expandafter{\romannumeral
2}}(22;6,6,6;5,2,2)$. By calculating, we know that
$\mathcal{E}(G_0)<\mathcal{E}(G_2)$.

Therefore, the proof is complete. \qed

From Theorem \ref{Thm3.8} and Lemmas \ref{lemma3.3},
\ref{lemma3.9} and \ref{lemma3.10}, we can easily
obtain the following theorem.
\begin{thm}\label{Thm3.11}
For any graph $G\in \Theta_{\uppercase\expandafter{\romannumeral
1}}(n;a,b,k;l_{1},l_{2};l_c)$ and $G\notin\Gamma_{1}$, there
exists a graph $H\in \Theta_{\uppercase\expandafter{\romannumeral
2}}(n;6,6,6;n-17,2,2)$ such that $G\preceq H$, and the equality holds
if and only if $G\cong H$.
\end{thm}

From Theorems \ref{Thm3.7} and \ref{Thm3.11}, we can obtain our main
result Theorem \ref{mainthm}.

\end{document}